\documentclass[11pt]{amsart}
\usepackage{amssymb,amsfonts,epsfig,subfigure}
\title[Index and Transversality]
  {Perturbed closed geodesics are periodic orbits: Index and Transversality}
\author{Joa Weber}
\address{Departement Mathematik, ETH Zentrum, 8092 Z\"urich, Switzerland}
\email{joa@math.ethz.ch, http://www.math.sunysb.edu/$^\sim$joa}
\date{\today}
\newtheorem{thm}{Theorem}[section]
\newtheorem{lem}[thm]{Lemma}
\newtheorem{prop}[thm]{Proposition}

\theoremstyle{definition}
\newtheorem{defn}[thm]{Definition}
\newtheorem{rmk}[thm]{Remark}

\newcommand{\1}{{{\mathchoice {\rm 1\mskip-4mu l} {\rm 1\mskip-4mu l}
  {\rm 1\mskip-4.5mu l} {\rm 1\mskip-5mu l}}}}
\newcommand{\comp}{{\scriptstyle \circ}}
\newcommand{\eps}{\epsilon}

\newcommand{\n}{\nabla}
\newcommand{\p}{\partial}

\newcommand{\smb}[1]{{\scriptstyle (} #1 {\scriptstyle )}}
\newcommand{\te}{\textstyle}
\newcommand{\AMC}{\MC{A}}
\newcommand{\BMC}{\MC{B}}

\newcommand{\CMC}{\MC{C}}
\newcommand{\CI}{C^\infty}

\newcommand{\EMC}{\MC{E}}
\newcommand{\FMC}{\MC{F}}

\newcommand{\LMC}{\MC{L}}

\newcommand{\MC}[1]{{\mathcal #1}}
\newcommand{\N}{{\mathbb N}}
\newcommand{\PMC}{\MC{P}}
\newcommand{\Q}{{\mathbb Q}}
\newcommand{\R}{{\mathbb R}}
\newcommand{\SMC}{\MC{S}}

\newcommand{\VMC}{\MC{V}}
\newcommand{\Z}{{\mathbb Z}}

\begin{document}

\begin{abstract}
  We study the classical action functional $\SMC_V$ on the free loop space
  of a closed, finite dimensional Riemannian manifold $M$ and the symplectic action
  $\AMC_V$ on the free loop space of its cotangent bundle. The critical points
  of both functionals can be identified with the set of perturbed closed geodesics in $M$. 
  The potential $V\in C^\infty(M\times S^1,\R)$
  serves as perturbation and we show that both functionals are Morse
  for generic $V$. In this case we prove that the Morse index of a
  critical point $x$ of $\SMC_V$ equals minus its Conley-Zehnder index
  when viewed as a critical point of $\AMC_V$ and if $x^*TM\to S^1$ is trivial.
  Otherwise a correction term $+1$ appears.
\end{abstract}


\maketitle
\tableofcontents
 
\section{Introduction and main results}

  We consider a smooth, compact Riemannian manifold $M$
of dimension $n$ and without boundary. The inner
product on the tangent space $T_xM$
is denoted by $\langle \cdot , \cdot \rangle$ 
and $g(x) : T_xM \to T_x^*M$ is the induced isomorphism.
We study the set of critical points of the classical action functional $\SMC_V$
on the space $\LMC M$ of free smooth loops in $M$ which is defined by
\begin{equation*}
  \SMC_V (x) :=
  \int_0^1 L(t,x(t), \dot x (t)) \: dt
\end{equation*}
The Lagrangian function $L=L_V : S^1 \times TM \to \R$ has the special form
kinetic minus potential energy 
$$
  L(t,x(t),\dot x(t)):=\frac{1}{2} \left| \dot x(t) \right|^2 -V(t,x(t))
$$
with a time periodic potential
$V \in C^\infty (S^1 \times M, \R)$.
Here and throughout we identify $S^1$ with $\R / \Z$ and represent the loop $x \in \LMC M$
as a smooth periodic function $x:\R \to M$ satisfying $x(t+1)=x(t)$.
The $L^2$-gradient of $\SMC_V$ is easily computed to be
\begin{equation*}
  grad \: \SMC_V (x)= -\n_t \dot x -\n V(t,x)
\end{equation*}
where $\n V$ denotes the gradient of $V$ with respect to the $x$-variable and where
$\n_t$ denotes covariant differentiation in direction $\dot x$
with respect to the Levi-Civita connection $\n$. The set of critical points of $\SMC_V$
is abbreviated by
\begin{equation} \label{eq:Crit}
  Crit = Crit \: \SMC_V= \{ x \in C^\infty(S^1,M) \mid -\n_t \dot x -\n V(t,x)=0 \}
\end{equation}
We call its elements \emph{perturbed closed geodesics}, since
it coincides with the set of closed geodesics on $M$ 
in the special case of constant potential $V$.

  On the other hand the symplectic action functional $\AMC_V$
on the free loop space $\LMC T^*M$ of $T^*M$ is defined as
\begin{equation*}
   \AMC_V (z)
   :=\int_0^1 \left( \langle y(t) , \dot x(t) \rangle - H(t,z(t)) \right) dt
\end{equation*}
where $z=(x,y)$ and the time dependent Hamiltonian function $H=H_V : S^1 \times T^*M \to \R$ is 
the Legendre transform of $L$; namely kinetic plus potential energy
\begin{equation} \label{eq:hamiltonian}
  H(t,x,y)=\frac{1}{2} \left| y \right|^2 + V(t,x) , \qquad (x,y) \in T^*M
\end{equation}
The set $Crit \: \AMC_V$ of critical points of $\AMC_V$ can be naturally identified 
with the set $Crit \: \SMC_V$
via the bijection
\begin{equation*}
  Crit \: \SMC_V \to Crit \: \AMC_V :
  x \mapsto z_x=(x,g(x) \dot x)
\end{equation*}
The inverse map is given by projecting onto $M$.
Note that $\SMC_V(x)=\AMC_V(z_x)$ for $x \in Crit$.

  Our first result asserts that the functionals
$\SMC_V$ and $\AMC_V$ are Morse for generic $V$. More precisely, in subsection \ref{subsec:hessians}
it will turn out that
the Hessian of $\SMC_V$ at a critical point $x$
can be represented by the perturbed Jacobi operator
$A_x$ in $L^2(S^1,x^*TM)$ with dense domain $W^{2,2}(S^1,x^*TM)$
given by
\begin{equation} \label{eq:A_x}
  A_x \: \xi = -\n_t \n_t \xi - R(\xi ,\dot x) \dot x -\n_\xi \n V_t(x)
\end{equation}
and where $R$ denotes the Riemann curvature tensor.
This operator is injective if and only if it is surjective.
A critical point $x$ of $\SMC_V$ is called \emph{nondegenerate} if
$A_x$ is injective. A function with nondegenerate
critical points only is called a \emph{Morse function}.
For $a \in \R$ let
$$
  \LMC^a M := \LMC M \cap \SMC_V^{-1}(-\infty,a] , \qquad
  Crit^a := Crit \cap \LMC^a M
$$
The space $\LMC^a T^*M$ is defined similarly with $\SMC_V$ replaced by $\AMC_V$.
Define the set $\VMC_{reg}$ of \emph{regular potentials} to be the set
of all $V\in \CI (S^1\times M,\R)$ such that $\SMC_V : \LMC M \to \R$
is Morse. Restricting $\SMC_V$ to $\LMC^a M$ we obtain the set
$$
  \VMC^a_{reg} := \{ V\in \CI (S^1\times M,\R) \mid
  A_x \; \mbox{is surjective} \; \forall x \in Crit^a \}
$$
In subsection \ref{subsec:hessians} we shall see that a critical point $x$
of $\SMC_V$ is nondegenerate if and only if the corresponding critical
point $z_x$ of $\AMC_V$ is. In other words $\SMC_V$ is Morse if and only if
$\AMC_V$ is.

\begin{thm} \label{thm:transversality}
  {\rm ({\bf Transversality})}
  i) For every $a \in \R$ the subset $\VMC^a_{reg} \subset \CI(S^1\times M,\R)$ 
  of regular potentials is open and dense.
  \newline
  ii) The subset of regular potentials $\VMC_{reg}\subset \CI(S^1\times M,\R)$
  is residual and therefore dense.
\end{thm}

  The set $\CI(S^1\times M,\R)$ is a complete metric space with the distance function
$$
  d(V_1,V_2) := \sum_{k=0}^\infty {\textstyle \frac{1}{2^k} \:
  \frac{\| V_1-V_2 \|_{C^k}}{1+\| V_1-V_2 \|_{C^k}} }
$$
A \emph{residual set} is by definition one which contains 
a countable intersection of open and dense sets.
According to Baire's category theorem \cite[section III.5]{RS80} such a set is dense.

  Choosing a regular potential $V$, we can assign two integers
to any $x \in Crit$ as follows: On the one hand the perturbed closed geodesic $x$
has a \emph{Morse index} $Ind(x)$, namely the number of negative eigenvalues of $A_x$
counted with multiplicities, and on the other hand
it is possible to interpret $z_x$ as a periodic orbit of the Hamiltonian
system on $T^*M$ with Hamiltonian (\ref{eq:hamiltonian}) 
and therefore define its Conley-Zehnder index $\mu_{CZ}(z_x)$. 
Since there exists a global Lagrangian splitting of $T(T^*M)$, the latter is well defined,
at least if $x^*TM \to S^1$ is trivial. Otherwise there is some choice involved
which leads to well definedness modulo $2$ only.
Our second result relates both indices.

\begin{thm} \label{thm:index}
  {\rm ({\bf Index})}
  Let $x \in Crit$ be nondegenerate, then
  $$
    \mu_{CZ} (z_x) = - Ind(x)
  $$
  if the bundle $x^*TM \to S^1$ is trivial and
  $$
    \mu_{CZ} (z_x) = - Ind(x) + 1
  $$
  if the bundle $x^*TM \to S^1$ is nontrivial.
\end{thm}

The relation between the Maslov index and the Morse index
of a closed geodesic has been studied first, as far as I know,
by Duistermaat \cite{Du76}. In the case of a closed geodesic
on a flat torus theorem \ref{thm:index} had been obtained
by Claude Viterbo \cite{Vi90} with a slightly different
definition of the Conley-Zehnder index (apart from different
normalizations):
due to the degeneracy of his action functional he considered
the Conley-Zehnder index of the linearized Hamiltonian flow
on the energy surface restricted to directions normal
to the trajectory.
\newline
A new feature arising in the present context is that $M$ is not required to be
orientable and so $x^*TM \to S^1$ may not be trivial. On the
other hand ${z_x}^*TT^*M \to S^1$ always is. We overcome this problem by
trivializing $x^*TM$ over $[0,1]$ and then apply to the induced trivialization
of ${z_x}^*TT^*M \to [0,1]$ an artificial half rotation to close up the frame
and obtain a trivialization over $S^1$.
Our choice of rotation is reflected in the formula by the term $\sigma(x)$.
Other choices lead to other odd integer multiples of this term.

\vspace{.3cm}
\noindent
{\bf Application.} Our results will be applied in \cite{We99}, \cite{SW01} and \cite{We01} 
to construct algebraic chain groups in the following way: Fix $a \in \R$, a regular 
potential $V$, and define
\begin{equation} \label{eq:sum}
  C_j^a (M,g,V;\Z_2) := \bigoplus_{\stackrel{\scriptstyle x \in Crit^a}{Ind(x)=j}} \
  Z_2 \langle x \rangle
\end{equation}
Finiteness of the sum is a consequence of transversality
combined with compactness as explained in appendix \ref{appsec:finite sum}.
In case the negative gradient flows of $\AMC_V$ and
$\SMC_V$ are Morse-Smale, we can count flow lines
between critical points of index difference $1$
to obtain chain maps $\p^F$ and $\p^M$, respectively. 
These flow lines are solutions of Floer's elliptic PDE on $T^*M$
and of the parabolic $L^2$-heat flow in the loop space of $M$, respectively.
It is Floer's theorem \cite[Thm. 4]{Fl89} that $\p^F \comp \p^F=0$; 
up to the additional argument
by Cieliebak \cite[theorem 5.4]{Ci94} 
in order to deal with noncompactness of $T^*M$.
In a forthcoming paper \cite{We01} it will be shown that $\p^M$
is a boundary operator too, whose homology represents the singular
homology of $\LMC^a M$.
In a joint research project with D. Salamon \cite{SW01} we will show that
both homology theories are naturally isomorphic
\begin{equation*} 
  HF_*^a (T^*M,H,J_g;\Z_2) \simeq
  HM_* (\SMC_{V,g},\LMC^a M;\Z_2)
\end{equation*}
As a consequence, the Floer homology of a cotangent bundle with a quadratic
Hamiltonian of the form (\ref{eq:hamiltonian}) is isomorphic to
the singular homology of the loop space. This result has been
obtained by Viterbo with different methods in \cite{Vi96}.
Our idea of proof is to obtain the heat equation
as an adiabatic limit of Floer's elliptic PDE.
The index theorem shows that the Morse index serves as a natural 
grading of Floer homology.

\begin{rmk} \label{rmk:sign-conventions}
  ({\bf Sign conventions and normalizations})
  The canonical symplectic structure $\omega_{can}$ on $T^*M$ is with respect to 
  natural coordinates $(x^1, \ldots ,x^n,
  y_1, \dots , y_n)$ given by $\sum_{i=1}^n dx^i \wedge dy_i$.
  The standard complex structure $J_0$ on $\R^{2n}$ and
  the signature of a symmetric matrix $S$ are defined by
  \[
    J_0:=\begin{pmatrix} 0&-\1 \\ \1&0 \end{pmatrix} , \qquad
    sign \: S := n^+(S) - n^-(S)
  \]
  where $n^{+(-)} (S)$ is the number of positive (negative)
  eigenvalues of $S$.
  The Conley-Zehnder index
  and the upward spectral flow are normalized as follows
  \begin{equation*}
    \begin{gathered}
      \mu_{CZ} ( t \mapsto e^{-t J_0 S})
      = -{\te \frac{1}{2}} \: sign \: S , \qquad
      \mu_{Spec} ( s \mapsto \arctan s ) =1
    \end{gathered}
  \end{equation*}
  where $t \in [0,1]$, $s \in \R$ and the symmetric matrix $S$ satisfies $\| S \|<2\pi$.
\end{rmk}

\noindent
{\bf Acknowledgements.} The author is most grateful to Kai Cieliebak, Dietmar Salamon
and Eduard Zehnder for stimulating discussions and comments.

\section{The index theorem} \label{sec:indexthm}

\subsection{The Hessians} \label{subsec:hessians}

  The critical points $x$ of the classical action $\SMC_V$ are the solutions of
the nonlinear equation
\begin{equation} \label{eq:CritS_V}
  -\n_t \p_t x - \n V_t(x) =0
\end{equation}
which follows from the formula
\begin{equation*}
  d \SMC_V(x) \xi = \int_0^1 \langle -\n_t \dot x - \n V_t(x) , \xi \rangle \: dt
  , \qquad \forall \xi \in C^\infty(S^1,x^*TM)
\end{equation*}
Linearization at a critical point $x$ leads to
the Hessian of $\SMC_V$ at $x$ given by
$$
  d^2 \SMC_V(x) (\xi,\eta) 
  = \langle A_x \xi , \eta \rangle_{L^2}
  =\int_0^1 \langle -\n_t \n_t \xi -R(\xi ,\dot x) \dot x
  -\n_\xi \n V_t(x) , \eta \rangle \: dt
$$

\noindent
In order to obtain this formula
introduce local coordinates and a smooth variation $x_\tau =(x_\tau^1, \ldots,
x_\tau^n)$ of $x$; i.e. $x_\tau$ depends smoothly on $\tau \in (-\delta,\delta)$,
$\delta >0$ small, and $x_0=x$, 
$\left. \frac{d}{d\tau} \right|_{\tau=0} x_\tau =:(\xi^1, \ldots ,\xi^n)$.
In (\ref{eq:CritS_V}) replace $x$ by $x_\tau$, apply 
$\left. \frac{d}{d\tau} \right|_{\tau=0}$ and
use the nonlinear equation as well as the representation of $R$ in local
coordinates via the Christoffel symbols of $\n$ to obtain (\ref{eq:A_x}).

  On the other hand the critical points of the symplectic action
$\AMC_V$ correspond precisely to the 1-periodic solutions
of the Hamiltonian system $(T^*M, \omega_{can}=-d\lambda, H_V)$.
Here $\lambda$ is the Liouville form which, in natural coordinates
$(x^1,\dots,x^n,y_1,\ldots,y_n)$ of $T^*M$  introduced below, 
is given by $\sum_{i=1}^n y_i dx^i$.
The Hamiltonian differential equation $\dot z (t) = X_{H_V}(t,z(t))$
can be expressed in the form
\begin{equation} \label{eq:CritA_V}
  \begin{split}
    \begin{pmatrix} \p_t x \\ \n_t y \end{pmatrix}
    =\begin{pmatrix} g(x)^{-1} y \\ -g(x) \n V_t(x) \end{pmatrix}
  \end{split}
\end{equation}
where $z(t)=(x(t),y(t))$ with $y(t)\in T_{x(t)}^*M$.
Note that here and throughout we identify $T_{z(t)}T^*M$ with
$T_{x(t)}M \oplus T_{x(t)}^*M$ by the isomorphism which takes
$\p_t z(t)$ to $(\p_t x(t), \n_t y(t))$. The equivalence of
(\ref{eq:CritS_V}) and (\ref{eq:CritA_V}) is then obvious.
\noindent
In order to linearize (\ref{eq:CritA_V}) at a zero $(x,y)$ introduce local
coordinates $(x^1,\dots,x^n)$ for $x$ on $M$ and obtain natural
coordinates for $y$ determined by $y=\sum_j y_j dx^j$.
Choose a smooth variation $(x_\tau^1 , \ldots , y^\tau_n)$ of $(x,y)$
and denote
$$
  (\xi^1,\dots,\xi^n,\tilde{\eta}_1,\ldots,\tilde{\eta}_n)
  :=\left. \frac{d}{d\tau} \right|_{\tau=0}
  (\xi_\tau^1,\dots,\xi_\tau^n,\tilde{\eta}_1^\tau,\ldots,\tilde{\eta}_n^\tau)
$$
It turns out that the $\xi^k$ transform as components of a vector $\xi \in T_xM$
under coordinate changes, but the $\tilde{\eta}_\ell$ do not have
any global meaning. However, the following combination with
the Christoffel symbols of the Levi-Civita connection $\n$
$$
  \eta_\ell:=\tilde{\eta}_\ell - \Gamma_{i\ell}^k(x) \xi^i y_k , \qquad \ell=1,\dots,n
$$
represents a covector $\eta=\sum_\ell \eta_\ell dx^\ell \in T_x^*M$.
Recall that $y_k$ denotes the $k^{th}$ component 
of the fibre part of the chosen zero $(x,y)$ of (\ref{eq:CritA_V})
with respect to natural coordinates.
Now linearization of (\ref{eq:CritA_V}) at a solution $z_x=(x,g(x) \p_t x)$
of (\ref{eq:CritA_V}) leads to the selfadjoint operator
$A_{z_x}$ in $L^2(S^1,x^*TM\oplus x^*T^*M)$ with dense domain $W^{1,2}(S^1,x^*TM \oplus x^*T^*M)$
which represents the Hessian of $\AMC_V$ at $z_x$ and is given by
\begin{equation} \label{eq:d^2A_V(z_x)}
  A_{z_x} \begin{pmatrix} \xi \\ \eta \end{pmatrix}
  =\begin{pmatrix} -g^{-1}(x) \n_t \eta -R(\xi,\dot x ) \dot x -\n_\xi \n V_t (x) \\
  g(x) \n_t \xi - \eta \end{pmatrix}
\end{equation}
This operator is clearly injective if and only if $A_x$ is:
a short calculation shows that $(\xi, g(x) \n_t \xi) \in ker \: A_{z_x}$ 
if and only if $\xi \in ker \: A_x$
and this proves

\begin{lem} \label{lem:AMC_V-Morse-SMC_V}
   $\AMC_V$ is a Morse function if and only if $\SMC_V$ is.
\end{lem}

\subsection{Orthonormal and unitary trivializations} \label{subsec:trivializations}

  In order to compare the Morse index of a critical point $x$ of $\SMC_V$ 
and the Conley-Zehnder index of the critical point $z_x$ of $\AMC_V$ it will be convenient to
view $x^*TM$ as vector bundle over $[0,1]$, choose an orthonormal trivialization
$\phi$ satisfying a certain boundary condition
and express the perturbed Jacobi operator $A_x$ with respect to this orthonormal frame.
This same frame $\phi$ will then be used to construct a unitary trivialization $\Phi_U$
for the trivial bundle $x^*TM \oplus x^*T^*M \to S^1$. In case $x^*TM \to S^1$ is nontrivial
some additional half rotation $U$ has to be applied in order to obtain indeed a frame over
$S^1$ rather than $[0,1]$.
More precisely, let
$$
  \sigma =\sigma(x):=
  \begin{cases}
    0 & \text{if $x^*TM$ is trivial}, \\
    1 & \text{otherwise}, \\
  \end{cases}
$$
and
$$
  E_\sigma := diag_n \left( (-1)^\sigma ,1,\ldots,1 \right) \in \R^{n \times n}
$$
Since $O(n)$ has precisely two connected components
one of which contains $E_0$ and the other one $E_1$, we may choose an orthonormal trivialization
\begin{equation} \label{eq:phi}
  \phi=\phi_\sigma : [0,1] \times \R^n \to x^*TM
\end{equation}
such that $\phi(1)=\phi(0) E_\sigma$.
Let $\dot \xi$ denote $\p_t \xi$, define the space
$$
  C_\sigma^\infty ( [0,1],\R^n )
  := \{ \xi \in C^\infty ( [0,1],\R^n ) \mid
  \xi(1)=E_\sigma \xi(0) , \; \dot \xi(1) = E_\sigma \dot \xi(0) \}
$$
and let $W^{2,2}_\sigma =W^{2,2}_\sigma ( [0,1],\R^n )$ and 
$L^2_\sigma=L^2_\sigma ( [0,1],\R^n )$ be
the closure of $C_\sigma^\infty ( [0,1],\R^n )$
with respect to the Sobolev  $W^{2,2}$- and $L^2$-norms, respectively.

  We define the model operator
$A^0 : L^2_\sigma \supset W^{2,2}_\sigma
\to  L^2_\sigma $
for $A_x$ with respect to the isometry induced by $\phi$ by
\begin{equation} \label{eq:A0}
    A^0 \xi 
    := \phi^{-1} A_x \: \phi \xi
    = -\p_t \p_t \xi - B \xi - Q \xi
\end{equation}
Induced covariant differentiation in direction $\dot x$
in the trivial bundle $[0,1]\times \R^n$
is of the form
\begin{equation} \label{eq:P}
  \phi^{-1}(t) \n_t \phi(t) \xi(t)= \p_t \xi(t) + P(t) \xi(t)
\end{equation}
where the connection potential $P$ is a family of skewsymmetric matrices.
Then the $L^2$-symmetric first order operator $B$ and the
family of symmetric matrices $Q$ in (\ref{eq:A0}) are given by
\begin{equation*}
  \begin{gathered}
    B \xi=(\p_t +P)^2 \xi -\p_t \p_t \xi
    =2P\p_t \xi + (\p_t P) \xi + P^2 \xi \\
    Q \xi = \phi^{-1} \bigl( R(\phi \xi , \dot x) \dot x
    +\n_{\phi \xi} \n V_t(x) \bigr)
  \end{gathered}
\end{equation*}

  For $t\in \R$ let now $U(t)$ be the $2n \times 2n$ matrix which represents
rotation of the $(1,n+1)$-coordinate plane
by the angle $\pi t$ and is the identity on all other coordinates
\begin{equation} \label{eq:U}
  U(t)
  =\begin{pmatrix} diag_n(\cos \pi t ,1,\ldots,1) & diag_n(-\sin \pi t ,0,\ldots,0) \\
  diag_n(\sin \pi t ,0,\ldots,0) & diag_n(\cos \pi t ,1,\ldots,1) \end{pmatrix}
\end{equation}
The orthonormal trivialization $\phi=\phi_\sigma$ introduced above
leads to a unitary trivialization of $x^*TM\oplus x^*T^*M$ over $[0,1]$ only, namely
$$
  \Phi =\Phi_\sigma := \begin{pmatrix} \phi_\sigma & 0 \\ 0 & {\phi_\sigma^*}^{-1} \end{pmatrix}
$$
Multiplication by $U(t) \in Sp(2n) \cap O(2n)$
then gives rise to a unitary trivialization of $x^*TM\oplus x^*T^*M$ over $S^1$
\begin{equation} \label{eq:PhiU}
  \Phi_U = \Phi_{U,\sigma} := \Phi_\sigma U^{-\sigma}
  \; : \; S^1 \times \R^{2n} \to x^*TM \oplus x^*T^*M
\end{equation}
In case $\sigma(x)=0$, $U^{-\sigma} =\1$ takes no effect. This is fine
since we already have $\Phi(1)=\Phi(0)$, because of $\phi(1)=\phi(0)$. 
Otherwise multiplication by $U^{-1}$ serves to close up the frame at $t=1$.
We may assume without loss of generality that the frame closes up smoothly:
If not, we modify $\phi$ such that $\n_t \phi=0$ near the ends of  $[0,1]$.
Of course we could choose in (\ref{eq:PhiU}) instead of $U$ some other power
$U^{2k+1}$, where $k \in \Z$ determines direction and
multiplicity of rotation. This would lead to a change of the Conley-Zehnder index
in definition \ref{def:CZindex} by $2k\sigma(x)$ as stated in
lemma \ref{lem:depU} $(ii)$.

  Finally, the operator $A_{z_x}$ is represented with respect to the unitary
frame $\Phi_U$ by the linear operator $A^1$ in $L^2(S^1,\R^{2n})$
with dense domain $W^{1,2}(S^1,\R^{2n})$. Calculation leads to
\begin{equation} \label{eq:A1}
  A^1 
  := {\Phi_U}^{-1} \: A_{z_x} \: \Phi_U
  = J_0 \p_t - S_U
\end{equation}
where
$$
  S_U = U^\sigma SU^{-\sigma} - J_0 U^\sigma \p_t (U^{-\sigma}) , \qquad
  S=\begin{pmatrix} Q & P \\ -P & \1 \end{pmatrix}
$$
and $Q,P$ are the families of matrices in (\ref{eq:A0}) and (\ref{eq:P}), respectively.
It is easy to check that $S$ is a family of symmetric matrices and
$S_U$ is symmetric with respect to the $L^2$-inner product.

\subsection{Morse index} \label{subsec:Morse}

  Let $x$ be a critical point of $\SMC_V$ and consider the perturbed Jacobi-operator
$A_x$ defined in (\ref{eq:A_x}) as an unbounded operator in 
$L^2=L^2(S^1,x^*TM)$
with dense domain $W^{2,2}=W^{2,2}(S^1,x^*TM)$. It is selfadjoint since it consists of 
the operator $\frac{d^2}{dt^2}$ on $S^1$ plus a bounded operator.

\begin{thm} \label{thm:morse}
  {\rm ({\bf Morse index theorem})}
  Let $x \in Crit$, then the Morse index $Ind (x)$
  and the nullity $Null (x) := dim \; ker \; A_x$ are finite. 
\end{thm}

\begin{proof}
Since $M$ is compact there exists a constant $C>0$
such that
\begin{equation*}
  \begin{split}
    \langle \xi , A_x \xi \rangle_{L^2}
   &=\| \n_t \xi \|_2^2 - \langle \xi , R(\xi,\dot x) \dot x \rangle_{L^2}
     -\langle \xi , \n_\xi \n V_t (x) \rangle_{L^2} \\
   &\ge \| \n_t \xi \|_2^2 - C \| \xi \|_2^2
  \end{split}
\end{equation*}
for all $\xi \in C^\infty (S^1,x^*TM)$.
For any $\rho > C$ it follows that the unbounded operator
$A_x+\rho$ in $L^2$ with dense domain $W^{2,2}$
is positive definite and so in particular injective. Moreover, it is
selfadjoint and therefore also surjective with real spectrum.
When viewed as a bounded operator from $W^{2,2}$ to $L^2$
the open mapping theorem guarantees existence of a bounded inverse.
Together with a standard compact Sobolev embedding 
we obtain that the resolvent operator is compact:
$$
  (A_x+\rho)^{-1} : L^2(S^1,x^*TM) \stackrel{\text{bd.}}{\longrightarrow} W^{2,2}(S^1,x^*TM)
  \stackrel{\text{cp.}}{\hookrightarrow} L^2(S^1,x^*TM)
$$
Compactness implies discrete spectrum, say $\{ 1/ \mu_j \}_{j\in \N}$, with
finite multiplicities and possible accumulation point at $0$.
We observe that $1/ \mu_j$ eigenvalue of $(A_x+\rho)^{-1}$
if and only if $\mu_j$ is an eigenvalue of $A_x+\rho$. We already know that
$\mu_j >0$ and conclude $\mu_j \to +\infty$ for $j\to \infty$.
Clearly,
$$
  (A_x+\rho)^{-1} \xi_j = {\textstyle \frac{1}{\mu_j}} \xi_j \quad 
  \Longleftrightarrow \quad 
  A_x \xi_j = (\mu_j -\rho) \xi_j
$$
Hence the eigenvalues $(\mu_j -\rho)$ of $A_x$ tend to $+\infty$
for $j\to \infty$ which proves the Morse index theorem.
\end{proof}

\subsection{Conley-Zehnder index} \label{subsec:CZ}

  Let us first illustrate the Conley-Zehnder index $\mu_{CZ}$ for a certain
class of paths in $Sp(2n)$, introduced in 1984 by Conley and Zehnder \cite{CZ84}. Later on we
shall give a precise definition of the more general Robbin-Salamon index $\mu_{RS}$ 
which will be more convenient to carry out calculations.

  Let $Sp^\pm(2n) = \{ Y \in Sp(2n) \mid \det (\1 -Y) \gtrless 0 \}$, $Sp^*(2n)=Sp^+(2n) \cup Sp^-(2n)$,
$\CMC(2n) = Sp(2n) \setminus Sp^*(2n)$ and $\SMC \PMC(2n)$ 
be the set of \emph{admissible paths}, which by definition means continous paths
$\gamma : [0,1] \to Sp(2n)$ such that $\gamma (0)=\1$ and $\gamma(1) \in Sp^*(2n)$.
The \emph{Maslov cycle} $\CMC(2n)$ is a codimension one algebraic subvariety of $Sp(2n)$ and it is
possible to interpret $\mu_{CZ} (\gamma)$ as the algebraic intersection number
of a generic path $\gamma \in \SMC \PMC(2n)$ with the Maslov cycle $\CMC(2n)$.
Generic means that $\gamma$ is of class $C^1$ and that the intersection of $\gamma$ with $\CMC(2n)$
is transversal for $t>0$.
Moreover, we need to assume that $\gamma$ departs from $\1$ at $t=0$ into $Sp^-(2n)$;
otherwise we homotop $\gamma$ within $\SMC \PMC(2n)$ to another path $\gamma^\prime$
satisfying the additional condition and define $\mu_{CZ}(\gamma)=\mu_{CZ}(\gamma^\prime)$.

\begin{figure}[ht]
  \centering
  \mbox{
    \subfigure[The path $\gamma_1$]
      {\epsfig{figure=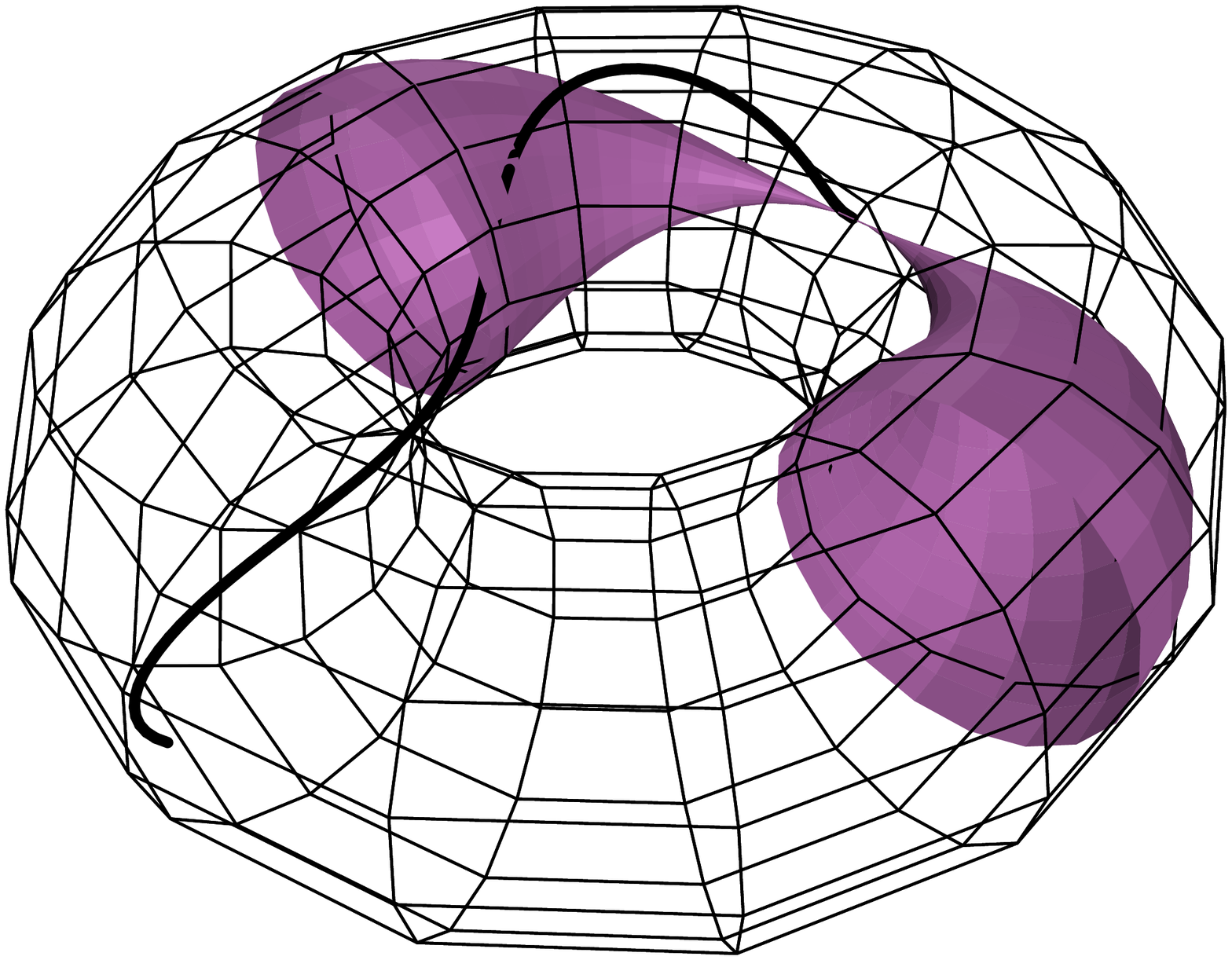,height=4cm,bbllx=18,bblly=180,bburx=594,bbury=630}} \quad
    \subfigure[The path $\gamma_2$]
      {\epsfig{figure=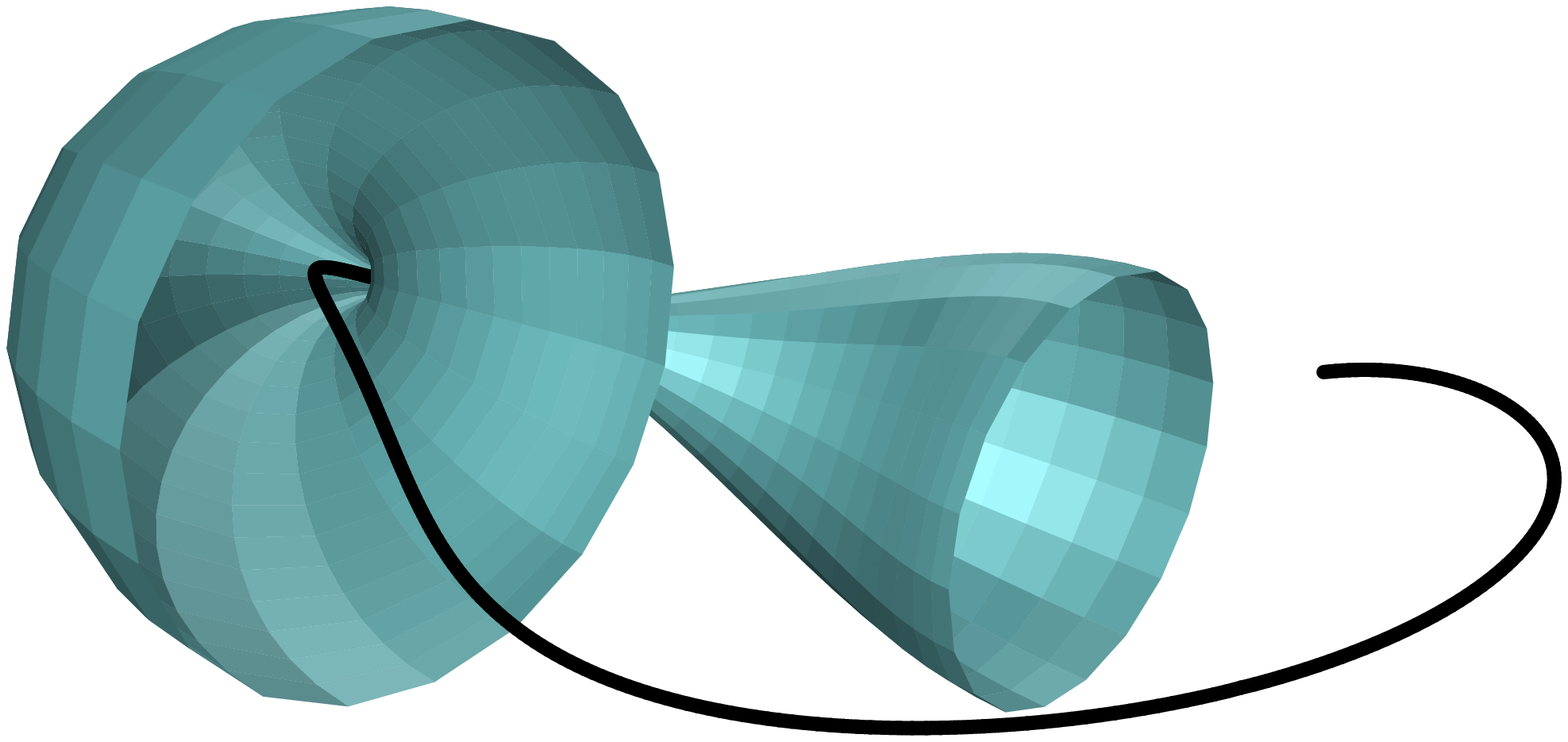,height=4cm,bbllx=18,bblly=180,bburx=594,bbury=630}}
  }
  \caption{Maslov cycle $\CMC(2)$ and paths with $\mu_{CZ}=+1$}
  \label{fig:mascyc}
\end{figure}

  For $n=1$ one can identify the symplectic linear group with the interior
of the solid $2$-torus \cite{GL58}; in this case $\CMC(2)$ has precisely
one singularity which corresponds to the identity matrix.
Figure \ref{fig:mascyc} (a) shows a numerical plot of the path 
$\gamma_1 : [0,1] \to Sp(2)$
$$
  \gamma_1(t)
  =\begin{pmatrix} \cos \pi t & - \sin \pi t \\ \sin \pi t & \cos \pi t \end{pmatrix}
  \begin{pmatrix} 1+t & 0 \\ 0 & (1+t)^{-1} \end{pmatrix}
$$
which has Conley-Zehnder index $+1$. In (b) the path $\gamma_2$ from (\ref{eq:gamma_2}) is shown
and it is important to notice that it departs from $\1$ at $t=0$
not into $Sp^-(2)$, but into $Sp^+(2)$. However, the path $\gamma_2$
is homotopic within $\SMC \PMC(2)$ to $\gamma_1$ and therefore its
Conley-Zehnder index is also $+1$.

  We observed in section \ref{subsec:hessians} that for any $x \in Crit$
we obtain a $1$-periodic orbit of the Hamiltonian system
$(T^*M,\omega_{can},H=H_V)$ by setting $z=z_x=(x,g(x) \p_t x)$ and every $1$-periodic orbit 
is obtained this way. Let $\varphi_t : T^*M \to T^*M$
denote the time-$t$-map generated by the Hamiltonian vector field
$X_H$, which is defined by $\omega_{can} (X_H ,\cdot ) =dH (\cdot )$, 
let $z_0=z(0)$ and $\Phi_U$ be the unitary trivialization
of $x^*TM \oplus x^*T^*M \to S^1$ introduced in (\ref{eq:PhiU}).
Then we obtain a smooth path of symplectic matrices by linearizing
the flow along the orbit $z(t)=\varphi_t z_0$
\begin{equation} \label{eq:linflow}
  \gamma (t) = \gamma_{x,U,\phi} (t) :=
  {\Phi_U(t)}^{-1} d\varphi_t (z_0) \Phi_U (0)
\end{equation}
Clearly $\gamma (0)=\1$ and the second condition $\gamma (1)
\in Sp^*$ is equivalent to $ker \; A_x =\{ 0 \}$, which reflects
our choice $V \in \VMC_{reg}$.

\begin{defn} \label{def:CZindex}
  For $x\in Crit$ such that $ker \: A_x = \{ 0 \}$ we define the Conley-Zehnder index
  of the $1$-periodic orbit $z_x=g(x) \p_t x$ by
  $$
    \mu_{CZ} (z_x) := \mu_{CZ} (\gamma_{x,U,\phi})
  $$
  where $\gamma_{x,U,\phi}$ is as in (\ref{eq:linflow}).
\end{defn}

\noindent
The following lemma shows that this definition
is actually independent of the choice of $\phi$. Its proof will be
given in subsection \ref{subsec:RS-index}.

\begin{lem} \label{lem:depU}
  (i) $\mu_{CZ} (\gamma_{x,U,\phi}) = \mu_{CZ} (\gamma_{x,U,\tilde{\phi}})$ for any
  $\phi , \tilde{\phi}$ as in (\ref{eq:phi}). \\
  (ii) $\mu_{CZ}(\gamma_{x,U^{2k+1}}) = 2 k \sigma(x) + \mu_{CZ}(\gamma_{x,U})$ for any 
  $k\in \Z$.
\end{lem}

\begin{rmk} \label{rmk:CZ}
  (1) In the standard setting of Floer homology,  namely that of a closed
  symplectic manifold $(N,\omega)$, one only considers contractible $1$-periodic
  orbits $z$ and trivializes $z^*TM \to S^1$ by filling in a disk. However, the Conley-Zehnder
  index of the resulting path in $Sp(2n)$ might depend on the homotopy class of the disk.
  An ambiguity arises precisely in the presence of spheres $\iota : S^2 \hookrightarrow N$
  with nonvanishing first Chern class $c_1(\iota^* TN)$; cf. the exposition \cite{Sa99}.
  Here this cannot happen, since
  $c_1(\iota^* TW)=0$ for any closed submanifold $\iota : W \hookrightarrow T^*M$;
  e.g. see \cite[section B.1.7]{We99} for two different proofs.
  \\
  (2) The reason that in the present case the Conley-Zehnder
  index is well defined for any closed orbit $x$ with $\sigma(x)=0$, contractible or not,
  is the existence of the global Lagrangian splitting
  $T_zT^*M \simeq T_xM \oplus T_x^*M$, $z \in T_x^*M$.
  However, as mentioned earlier, a $(mod \: 2)$-ambiguity
  arises in case of nontrivial bundles $x^*TM \to S^1$.
\end{rmk}

  In order to prove the index theorem it will be useful to represent
the path $\gamma$ in (\ref{eq:linflow}) as solution $\Psi=\Psi_U$ 
of the initial value problem stated in the following lemma.

\begin{lem} \label{lem:IVP}
  The path of symplectic matrices $\gamma$ 
  defined in (\ref{eq:linflow}) equals the fundamental solution $\Psi$
  of the initial value problem
  $$
    \p_t \Psi = -J_0 S_U \Psi , \qquad \Psi (0) =\1
  $$
  where $S_U$ is the family of symmetric matrices from (\ref{eq:A1}).
\end{lem}

\begin{proof}
  To show equality of $\gamma$ and $\Psi$ we choose the following smooth variation of
  $z(t)=\varphi_t (z_0)$
  $$
    z^\tau
    =\varphi_t (exp_{x_0} \tau \xi_0 , y_0 +\tau \eta_0 ) , \qquad 
    (\xi_0 ,\eta_0) \in T_{x_0}M \oplus T_{x_0}^*M
  $$
  in order to linearize (\ref{eq:CritA_V}) as described in subsection \ref{subsec:hessians}.
  With this choice
  $z^0(t)=\varphi_t (z_0)=z(t)=(x(t),y(t))$ and
  $$
    \begin{pmatrix} \xi (t) \\ \eta (t) \end{pmatrix}
    := \left. {\textstyle \frac{d}{d \tau} } \right|_{\tau=0} z^\tau (t)
    =d\varphi_t (z_0) \begin{pmatrix} \xi_0 \\ \eta_0 \end{pmatrix}
    \in C^\infty ([0,1],x^*TM \oplus x^*T^*M)
  $$
  is a zero of the linear equations (\ref{eq:d^2A_V(z_x)}), which however
  might not close up at time $t=1$.
  In other words the linearized flow along a flow line (zero of the nonlinear equations)
  provides a zero of the linearized equations.
  With respect to the unitary frame $\Phi_U$ these are given by
  $$
    0=(J_0 \p_t -S_U) {\Phi_U}^{-1} \begin{pmatrix} \xi \\ \eta \end{pmatrix}
    =(J_0 \p_t -S_U) \gamma
  $$
  and this proves the claim.
\end{proof}

\begin{rmk} \label{rmk:linearized-flow}
  Since the flow applied to a point
  produces a path, the natural domain for the linearized equations
  along such a path actually is the space of vector fields along the path.
  If we restrict our attention to those paths which happen
  to be $1$-periodic, then the natural domain becomes
  the set of vector fields along loops and we obtain the operator 
  $A_{z_x}$ from (\ref{eq:d^2A_V(z_x)}).

    To summarize, we have that any $(\xi,\eta)$ in the kernel of
  the linearized flow along a flow line $z=(x,y)$
  is of the form $d\varphi_t (z(0)) (\xi_0,\eta_0)$
  where $(\xi_0,\eta_0) \in T_{x(0)}M \oplus T_{x(0)}^*M$.
  In the periodic case the vector fields are required to close up at time $t=1$
  and so the corresponding kernel is isomorphic to a subspace, namely
  $$
    ker \: A_{z_x}
    \simeq ker \: \bigl( \1-d\varphi_1(z(0)) \bigr)
  $$
\end{rmk}

\subsection{Robbin-Salamon index} \label{subsec:RS-index}

  In order to prove the index theorem \ref{thm:index} as well as
lemma \ref{lem:depU} we shall briefly recall another
Maslov-type index $\mu_{RS}$, called \emph{Robbin-Salamon index} \cite{RS93}.
It is defined for any continuous path of Lagrangian subspaces 
$\Lambda : [a,b] \ni t \mapsto \Lambda(t)$
of a given symplectic vector space $(V,\Omega)$
with respect to a fixed Lagrangian subspace $\Lambda_0$.

  For now let us assume the path $\Lambda$ is of class $C^1$.
Call $t_i \in [a,b]$ a \emph{crossing}
if $\Lambda (t_i) \cap \Lambda_0 \not= \{ 0\}$. 
For any such $t_i$ there is a quadratic form
on $\Lambda(t_i)$: Pick any Lagrangian complement $W$ of $\Lambda(t_i)$ and for
$v \in \Lambda(t_i)$ and sufficiently small $\eps>0$ define
$w(\eps) \in W$ by the condition $v+w(\eps) \in \Lambda(t_i+\eps)$. Then
$$
  \hat{Q}(\Lambda (t_i),\p_t \Lambda (t_i)) \: v 
  := {\te \left. \frac{d}{d\eps} \right|_{\eps=0}}
  \Omega (v,w (\eps))
$$
is a quadratic form on $\Lambda(t_i)$, which is independent of the choice of $W$
\cite[theorem 1.1]{RS93}.
The \emph{crossing form at $t_i$} is its restriction to $\Lambda(t_i) \cap \Lambda_0$ 
$$
  \Gamma(\Lambda,\Lambda_0,t_i) := \hat{Q}(\Lambda (t_i),\p_t \Lambda (t_i)) 
  \mid_{\Lambda(t_i) \cap \Lambda_0}
$$
and $t_i$ is called a \emph{regular crossing} if its crossing form is nonsingular.
A path $\Lambda$ is called \emph{regular} if all its crossings are regular.

\begin{defn} \label{def:RS}
  The Robbin-Salamon index of a regular Lagrangian path $\Lambda$ is defined to be
  $$
    \mu_{RS} (\Lambda,\Lambda_0)
    = \frac{1}{2} sign \: \Gamma(\Lambda,\Lambda_0,a)
    + \sum_{a<t<b} sign \: \Gamma(\Lambda,\Lambda_0,t) 
    + \frac{1}{2} sign \: \Gamma(\Lambda,\Lambda_0,b)
  $$
  where the sum runs over all crossings $t$.
\end{defn}

We remark that in \cite[section 2]{RS93} the following is shown: Any two regular 
Lagrangian paths which are homotopic with fixed
endpoints have the same Robbin-Salamon index and every continuous
Lagrangian path is homotopic with fixed endpoints to a regular one.
As a consequence one can define the Robbin-Salamon index for every continuous path.

  The most important property of $\mu_{RS}$ in this text will be its
\emph{catenation property}, which means that $\mu_{RS}$
is additive under composition of paths (with respect to decompositions of
the parameter domain $[a,b]$).
Moreover, for any path of symplectic matrices $\Psi \in \SMC \PMC(2n)$ the Robbin-Salamon
index reproduces the Conley-Zehnder index \cite[remark 5.3]{RS93} as follows
\begin{equation} \label{eq:RS=CZ}
  \mu_{RS} (Graph \: \Psi , \Delta) = \mu_{CZ} (\Psi)
\end{equation}
Here $\Delta$ denotes the diagonal in the symplectic vector space
$(\R^{2n \times 2n} , -\omega_0 \oplus \omega_0 )$
and $\omega_0 (\cdot , \cdot ) = \langle J_0 \cdot , \cdot \rangle$
is the standard symplectic form on $\R^{2n}$.
The \emph{loop property} of the Conley-Zehnder index \cite{DS94}
translates into 
\begin{equation} \label{eq:loopRS}
  \mu_{RS} ( Graph \: \Theta \Psi , \Delta)
  = \mu_{RS} ( Graph \: \Theta , \Delta)
  + \mu_{RS} ( Graph \: \Psi , \Delta)
\end{equation}
for any path $\Psi \in \SMC \PMC(2n)$ and any loop $\Theta : S^1 \to Sp(2n)$.

\begin{lem} \label{lem:loop-RS}
  For every unitary loop 
  $$
    \Theta= \begin{pmatrix} X & -Y \\ Y & X \end{pmatrix} : S^1 \to Sp(2n) \cap O(2n)
  $$
  it holds
  $$
    \mu_{RS} ( Graph \: \Theta , \Delta )
    = 2 \deg \Bigl[ \det (X+iY) : S^1 \to S^1 \Bigr]
  $$
\end{lem}

\noindent
We only sketch a proof: 
The idea is to shows that the right hand side satisfies those axioms
which determine the left hand side uniquely. These are the
direct sum, the normalization and the weak homotopy axioms.

  Finally let us derive a formula for the crossing form
in case $\Lambda (t)=Graph \: \Psi (t) \subset (\R^{2n \times 2n} , -\omega_0 \oplus \omega_0 )$,
$\Lambda_0 = \Delta$ and $\Psi : [a,b] \to Sp(2n)$ of class $C^1$.
First observe that $\Psi$ determines a path of symmetric matrices by
$$
  S(t) = J_0 (\p_t \Psi(t)) \Psi(t)^{-1}
$$
Let $t_i \in [a,b]$ be a crossing, $v \in Graph \: \Psi(t_i) \cap \Delta$ and
$W=0 \times \R^{2n}$. Then $v=(\zeta , \Psi(t_i) \zeta)=(\zeta , \zeta)$
for some $\zeta \in \R^{2n}$, $w(\eps)=(0,w_2(\eps))$
and the condition $v+w(\eps) \in Graph \: \Psi (t_i+\eps)$
leads to $w_2(\eps)=\Psi(t_i+\eps) \zeta - \zeta$. We calculate
\begin{equation*}
  \begin{split}
    \hat{Q}(\Lambda (t_i) , \p_t \Lambda (t_i)) \: v
   &={\te \left. \frac{d}{d\eps} \right|_{\eps=0}}
    (-\omega_0\oplus\omega_0) \Bigl( (\zeta,\zeta) , (0,w_2(\eps)) \Bigr) \\
   &={\te \left. \frac{d}{d\eps} \right|_{\eps=0}}
    \Bigl( -\omega_0 (\zeta,0) 
    + \omega_0 (\zeta,\Psi (t_i+\eps) \zeta -\zeta) \Bigr) \\
   &=\omega_0 (\zeta, \p_t \Psi (t_i) \zeta) \\
   &=-\langle \zeta ,S(t_i) \Psi(t_i) \zeta \rangle \\
   &=-\langle \zeta, S(t_i) \zeta \rangle ,
  \end{split}
\end{equation*}
and, identifying $(Graph \: \Psi(t_i)) \cap \Delta$ and $ker \: (\1 -\Psi (t_i))$
by $(\zeta,\zeta) \mapsto \zeta$, we obtain
\begin{equation} \label{eq:crossform=S}
  \Gamma(Graph \: \Psi,\Delta,t_i) \: (\zeta,\zeta) 
  =-\langle \zeta ,S(t_i) \zeta \rangle, \qquad \zeta \in ker \: \bigl( \1 -\Psi(t_i) \bigr)
\end{equation}

Now we are in position to prove lemma \ref{lem:depU}.

\begin{proof} ({\sc of Lemma} \ref{lem:depU})
  {\bf ad (i)} Consider the transition maps $\psi = \phi^{-1} \tilde{\phi} : [0,1] \to O(n)$
  and $\Psi = \Phi^{-1} \tilde{\Phi} : [0,1] \to Sp(2n) \cap O(2n)$, then
  \begin{equation*}
    \begin{split}
      \mu_{CZ} &(\gamma_{x,U,\phi} )
      =\mu_{CZ} ( U^\sigma \Psi U^{-\sigma} \gamma_{x,U,\tilde{\phi}}
      U(0)^\sigma \Psi(0)^{-1} U(0)^{-\sigma} ) \\
     &=\mu_{RS} ( Graph \: U^\sigma \Psi U^{-\sigma} , \Delta )
      +\mu_{CZ} (\gamma_{x,U,\tilde{\phi}} )
      +\mu_{RS} ( Graph \: \Psi(0)^{-1} , \Delta ) \\
     &=\mu_{CZ} (\gamma_{x,U,\tilde{\phi}} ) 
    \end{split}
  \end{equation*}
  where in the second equality we used (\ref{eq:RS=CZ}) as well as the
  loop property (\ref{eq:loopRS}). Since $\Psi(0)^{-1}$ is a constant loop
  the corresponding term vanishes. It remains to show
  $\mu_{RS} ( Graph \: U^\sigma \Psi U^{-\sigma} , \Delta )=0$:
  Since $U^\sigma \Psi U^{-\sigma} : [0,1] \to Sp(2n) \cap O(2n)$ 
  we can write
  $$
    U^\sigma \Psi U^{-\sigma} = \begin{pmatrix} X & -Y \\ Y & X \end{pmatrix}
    \; \text{where} \; X+iY \in U(n)
  $$
  and by lemma \ref{lem:loop-RS} the following is true:
  \begin{equation} \label{eq:RS-deg}
    \mu_{RS} ( Graph \: U^\sigma \Psi U^{-\sigma} , \Delta )
    = 2 \deg \Bigl[ \det (X+iY) : S^1 \to S^1 \Bigr]
  \end{equation}

  \noindent
  Now in case $\sigma=0$ we have $X=\psi, Y=0$ and are done. If
  $\sigma=1$ a calculation shows that
  $$
    \det \left( X+iY \right)
    =\det \begin{pmatrix}
      \psi_{11} & e^{-i\pi t} \psi_{12} & \hdots & e^{-i\pi t} \psi_{1n} \\
      e^{i\pi t} \psi_{21}  & \psi_{22} & \hdots & \psi_{2n} \\
      \vdots & \vdots & & \vdots \\
      e^{i\pi t} \psi_{n1}  & \psi_{n2} & \hdots & \psi_{nn}
     \end{pmatrix}
    = \det \psi
    = \pm 1
  $$
  {\bf ad (ii)}
  \begin{equation*}
    \begin{split}
      \mu_{CZ} &(\gamma_{x,U^{2k+1}} )
      =\mu_{CZ} ( U^{2k\sigma} \gamma_{x,U} U(0)^{-2k\sigma} ) \\
     &=\mu_{RS} ( Graph \: U^{2k\sigma} , \Delta )
      +\mu_{CZ} (\gamma_{x,U} )
      +\mu_{RS} ( Graph \: U(0)^{-2k\sigma} , \Delta ) \\
     &=2\sigma \deg 
      \left( t \mapsto e^{\pi i 2kt} \right)
      +\mu_{CZ} (\gamma_{x,U} ) \\
     &=2k\sigma + \mu_{CZ} (\gamma_{x,U} )
    \end{split}
  \end{equation*}
  where we used (\ref{eq:RS=CZ}) as well as the loop property (\ref{eq:loopRS})
  in the second equality. The third one follows from lemma \ref{lem:loop-RS}.
\end{proof}

\subsection{Spectral flow and Conley Zehnder index} \label{subsec:specflow-CZ}

  We provide the main tool to prove the index theorem, namely theorem \ref{thm:specflow} below,
which relates the spectral flow and the Conley-Zehnder index.

  Pick $T>0$ and smooth two-parameter families of matrices
$Q,P:[-T,T]\times[0,1] \to \R^{n \times n}$ such that
$$
  Q_\lambda(t)^T=Q_\lambda(t) , \qquad Q_\lambda(1)= E_\sigma Q_\lambda(0) E_\sigma^{-1}
$$
and
$$
  {P_\lambda(t)}^T=-P_\lambda(t) , \qquad P_\lambda(1)= E_\sigma P_\lambda(0) E_\sigma^{-1}
$$
For $\lambda \in [-T,T]$ consider the family $A_\lambda$ of selfadjoint operators
in $L^2_\sigma$ with dense domain $W^{2,2}_\sigma$ defined by
\begin{equation} \label{eq:Ilambda}
  A_\lambda ( \xi )
  := - \p_t \p_t \xi - B_\lambda \xi - Q_\lambda \xi
\end{equation}
where $L^2_\sigma$ and $W^{2,2}_\sigma$ 
have been introduced in section \ref{subsec:Morse} and
the $L^2$-symmetric family of first order operators $B_\lambda$ is defined by
\begin{equation} \label{eq:alambda}
  B_\lambda \xi
  =2P_\lambda \dot \xi + \dot P_\lambda \xi + {P_\lambda}^2 \xi
\end{equation}
The term $\p_t \p_t +B_\lambda$ represents
the second covariant derivative $\n_t \n_t$ with respect to an orthonormal frame.
To indicate the particular choices $\lambda=\mp T$ we frequently will use
the simpler notation $\mp$.

    Let us now define the upward spectral flow of the family
  $A_\lambda$, $\lambda \in [-T,T]$.
  Roughly speaking, it counts the number of eigenvalues
  changing sign from minus to plus during the deformation minus the ones
  changing sign in the opposite way.
  The real number $\lambda$ is called a \emph{crossing}, if $A_\lambda$ is not injective.
  In this case, following the exposition in \cite{RS95},
  we define the \emph{crossing operator}
  \begin{equation} \label{eq:crossingoperator-spec}
    \Gamma ( \{ A_\lambda \}_{\lambda \in [-T,T]} , \lambda)
    =P_\lambda^\perp (\p_\lambda A_\lambda) P_\lambda^\perp \mid_{Ker \: A_\lambda}
  \end{equation}
  where $P_\lambda^\perp : L^2_\sigma \to L^2_\sigma$ denotes
  the orthogonal projection onto $Ker \: A_\lambda$.
  We emphasize that, despite the similarity of notation, the object $P_\lambda^\perp$
  is entirely different from $P_\lambda$ and that $sign$ below denotes
  the signature of a quadratic form.
  A crossing is called \emph{regular}, if its crossing operator is nonsingular.
  The \emph{spectral flow} is characterized axiomatically in \cite{RS95}.
  In case all crossings are regular we may use lemma 4.27 in \cite{RS95}
  to actually \emph{define}
  $$
    \mu_{Spec}(\{ A_\lambda \}_{\lambda \in [-T,T]})
    =\sum_\lambda sign \: \Gamma (\{ A_\lambda \}_{\lambda \in [-T,T]} ,\lambda)
  $$
  Note that the sum is over all crossings $\lambda$ and that there are only finitely many
  of them in view of their regularity.

\noindent
In what follows we need to assume injectivity of $A_\mp$:
\begin{equation} \label{eq:injectivity-assumption}
  \mbox{\emph{assume that $Ker \: A_\mp =\{ 0 \}$}}
\end{equation}
Let us now state the main theorem of this section.

\begin{thm} \label{thm:specflow} 
  Let $\{ A_\lambda \}_{\lambda \in [-T,T]}$ be a regular family 
  as in (\ref{eq:Ilambda}) and assume injectivity of $A_\mp$.
  Let $\Psi_{\mp,U}$ be the symplectic paths associated to $A_\mp$
  by the equivalence of statements (S1) and (S4) below, then the upward spectral flow
  of the family is given by
  $$
    \mu_{Spec} ( \{ A_\lambda \}_{\lambda \in [-T,T]} )
    = \mu_{CZ} ( \Psi_{+,U} ) 
    - \mu_{CZ} ( \Psi_{-,U} )
  $$
\end{thm}

\noindent

  Equivalence of the subsequent four statements is fairly easy, but
nevertheless crucial in the proof of the theorem.
Together they show how the operator $A_\lambda$ leads to a 
path $\Psi_{\lambda,U}$ of symplectic matrices.
Throughout let $\xi \in C^\infty ([0,1],\R^n)$,
denote $\eta = \dot \xi + P_\lambda \xi$
and
$$
  S_\lambda = S_\lambda (t) 
  = \begin{pmatrix} Q_\lambda & P_\lambda \\ -P_\lambda & \1_n \end{pmatrix} , \qquad
  \begin{pmatrix} \xi_U \\ \eta_U \end{pmatrix}
  = U^\sigma \begin{pmatrix} \xi \\ \eta \end{pmatrix}
$$
where the rotation $U$ is defined in (\ref{eq:U}) and $\sigma \in \{ 0 , 1 \}$.

\vspace{.2cm}
\noindent
{\bf (S1)}
  $$
    -\p_t \p_t \xi - B_\lambda \xi - Q_\lambda \xi = 0
  $$

\vspace{.2cm}
\noindent
{\bf (S2)}
  $$
    \p_t \begin{pmatrix} \xi \\ \eta \end{pmatrix}
    =-J_0 S_\lambda \begin{pmatrix} \xi \\ \eta \end{pmatrix} , \qquad
  $$

\vspace{.2cm}
\noindent
{\bf (S3)}
  $$
    \begin{pmatrix} \xi (t) \\ \eta (t) \end{pmatrix}
    =\Psi_\lambda (t) \begin{pmatrix} \xi (0) \\ \eta (0) \end{pmatrix}
  $$
  and the fundamental solution $\Psi_\lambda (t)$ is determined by
  $$
    \p_t \Psi_\lambda = -J_0 S_\lambda \Psi_\lambda , \qquad
    \Psi_\lambda (0) = \1_{2n}
  $$

\vspace{.2cm}
\noindent
{\bf (S4)}
  $$
    \begin{pmatrix} \xi_U \\ \eta_U \end{pmatrix} (t)
    =\Psi_{\lambda,U} (t) \begin{pmatrix} \xi_U \\ \eta_U \end{pmatrix} (0)
  $$
  and the fundamental solution $\Psi_{\lambda,U} (t)$ is determined by
  $$
    \p_t \Psi_{\lambda,U} = -J_0 S_{\lambda,U} \Psi_{\lambda,U} , \qquad
    \Psi_{\lambda,U} (0) = \1_{2n}
  $$
  $$
    S_{\lambda,U} = U^\sigma S_\lambda U^{-\sigma} - J_0 U^\sigma \p_t (U^{-\sigma})
  $$

\begin{rmk} \label{rmk:welldefCZ}
  If in addition we require $\xi$ in (S1) to satisfy the boundary conditions
  $\xi(1) = E_\sigma \xi(0)$ and $\dot \xi(1) = E_\sigma \dot \xi(0)$, 
  i.e. $\xi \in ker \: A_\lambda$, it turns out that in (S3)
  $$
    \begin{pmatrix} \xi (0) \\ \eta (0) \end{pmatrix}
    \in ker \: (E_\sigma \oplus E_\sigma - \Psi_\lambda(1) )
  $$
  and in (S4)
  $$
    \begin{pmatrix} \xi_U (0) \\ \eta_U (0) \end{pmatrix}
    \in ker \: (\1 - \Psi_{\lambda,U} (1) )
  $$
  In view of the injectivity assumption (\ref{eq:injectivity-assumption})
  the symplectic paths $\Psi_{\mp,U}$ end outside the Maslov
  cycle and therefore are elements of the 
  set $\SMC \PMC(2n)$ for which the Conley-Zehnder index is defined.
\end{rmk}

\begin{proof} ({\sc of Theorem \ref{thm:specflow}})
  For each $\lambda \in [-T,T]$ there is the path of symplectic matrices 
  $\Psi_{\lambda,U}$ defined in (S4). Together these give rise to a 2-parameter family
  of Lagrangian subspaces
  $$
    \Lambda_\lambda (t) = Graph \: \Psi_{\lambda,U} (t)
    \subset (\R^{2n} \times \R^{2n}, -\omega_0 \oplus \omega_0 )
  $$
  and we consider the Robbin-Salamon index $\mu_{RS}(\Gamma,\Delta)$
  of the loop of Lagrangians $\Gamma$ along the boundary
  of the parameter domain $[0,1] \times [-T,T]$.
  \begin{figure}[ht]
    \centering\epsfig{figure=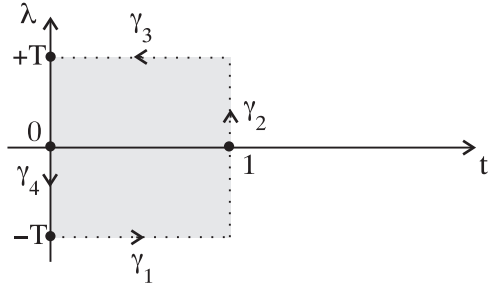}
    \caption{Contractible loop of Lagrangian subspaces}
    \label{fig:2parlagr}
  \end{figure}
  We may break down the loop $\Gamma$
  into the four subpaths $\gamma_1 ,\ldots, \gamma_4$
  indicated in figure \ref{fig:2parlagr}.
  $\Gamma$ is contractible and so its Robbin-Salamon index is zero.
  Together with the catenation property of $\mu_{RS}$ this leads to
  $$
    0=\mu_{RS} (\Gamma,\Delta)
    =\sum_{i=1}^4 \mu_{RS} (\gamma_i,\Delta)
  $$
  Since $\Psi_{\lambda,U} (0)=\1_{2n}$, $\gamma_4 \equiv \Delta$
  and so $\mu_{RS} (\gamma_4,\Delta)$=0.
  Moreover, as indicated in (\ref{eq:RS=CZ}) the Robbin-Salamon index of 
  $\gamma_1$ is precisely the Conley-Zehnder
  index of the corresponding symplectic path. We get a minus sign
  in the case of $\gamma_3$, because of its reversed orientation
  \begin{equation*}
    \begin{split}
      \mu_{RS}(\gamma_1,\Delta)
     &=\mu_{RS} (Graph \: \Psi_{-,U} ,\Delta)=\mu_{CZ}(\Psi_{-,U}) \\
       \mu_{RS}(\gamma_3,\Delta)
     &=-\mu_{RS} (Graph \: \Psi_{+,U} ,\Delta)=-\mu_{CZ}(\Psi_{+,U}) \\
    \end{split}
  \end{equation*}
  It remains to show
  \begin{equation} \label{eq:RS=specflow}
    \mu_{RS} (\gamma_2,\Delta)=\mu_{spec} (\{ A_\lambda \}_{\lambda \in [-T,T]} )
  \end{equation}

    In order to prove (\ref{eq:RS=specflow}) we need to derive two identities.
  Recall from (S4) that $\Psi_{\lambda,U}$ is determined by
  $$
    \p_t \Psi_{\lambda,U} (t)=-J_0 S_{\lambda,U}(t) \Psi_{\lambda,U} (t)
    , \qquad \Psi_{\lambda,U} (0) = \1
  $$
  Now for fixed $t \in [0,1]$ the path $\lambda \mapsto \Psi_{\lambda,U} (t)$,
  $\lambda \in [-T,T]$, leads to another path of symmetric matrices
  $$
    \hat S_t (\lambda) = J_0 (\p_\lambda \Psi_{\lambda,U}(t)) \Psi_{\lambda,U}(t)^{-1}
  $$
  Note that $\hat S_0(\lambda) \equiv 0$, since $\Psi_{\lambda,U}(0)\equiv \1$.
  Using these equations we obtain
  \begin{equation*}
    \begin{split}
      \p_t \bigl( \Psi_{\lambda,U}&(t)^T \: 
      \hat{S}_t(\lambda) \: \Psi_{\lambda,U}(t) \bigr) \\
     &=\bigl(\p_t \Psi_{\lambda,U}(t)^T \bigr) \: \hat{S}_t(\lambda) \: 
      \Psi_{\lambda,U}(t)
      +\Psi_{\lambda,U}(t)^T \: \p_t \bigl(J_0 \p_\lambda \Psi_{\lambda,U}(t) \bigr) \\
     &=\Psi_{\lambda,U}(t)^T \: S_{\lambda,U}(t) \: J_0 \:
      \hat{S}_t(\lambda) \: \Psi_{\lambda,U}(t) \\
     &\quad -\Psi_{\lambda,U}(t)^T \:  J_0 \: 
      \p_\lambda \bigl(J_0 \: S_{\lambda,U}(t) \: \Psi_{\lambda,U}(t) \bigr) \\
     &=-\Psi_{\lambda,U}(t)^T \: S_{\lambda,U}(t) \:
      \p_\lambda \Psi_{\lambda,U}(t)
      +\Psi_{\lambda,U}(t)^T \: \bigl(\p_\lambda S_{\lambda,U}(t) \bigr) \:
      \Psi_{\lambda,U}(t) \\
     &\quad +\Psi_{\lambda,U}(t)^T \: S_{\lambda,U}(t) \:
      \p_\lambda \Psi_{\lambda,U}(t) \\
     &=\Psi_{\lambda,U}(t)^T \: \bigl( \p_\lambda S_{\lambda,U}(t) \bigr) \:
      \Psi_{\lambda,U}(t)
    \end{split}
  \end{equation*}
  Integration over $t$ from $0$ to $1$, together with $\hat S_0(\lambda) \equiv 0$,
  leads to the first identity we are looking for, namely
  \begin{equation} \label{eq:identity1}
    \Psi_{\lambda,U}(1)^T \: \hat{S}_1(\lambda) \: 
    \Psi_{\lambda,U} (1)
    =\int_0^1 \Psi_{\lambda,U}(t)^T \: \bigl( \p_\lambda S_{\lambda,U}(t) \bigr) \:
    \Psi_{\lambda,U}(t) \: dt
  \end{equation}

  \noindent
  Use skew-symmetry of $P_\lambda$ and 
  $\p_\lambda P_\lambda$ as well as 
  $
    \eta = \dot \xi + P_\lambda \xi
  $
  and
  $$
    \p_\lambda B_\lambda
    =2(\p_\lambda P_\lambda) \p_t
    +\p_\lambda \dot P_\lambda +(\p_\lambda P_\lambda) P_\lambda
    +P_\lambda (\p_\lambda P_\lambda)
  $$
  to obtain the second identity
  \begin{equation} \label{eq:identity2}
    \begin{split}
      \int_0^1 
     &\left\langle \begin{pmatrix} \xi \\ \eta \end{pmatrix} ,
      \begin{pmatrix} -\p_\lambda B_\lambda & \p_\lambda P_\lambda \\
      -\p_\lambda P_\lambda & 0 \end{pmatrix}
      \begin{pmatrix} \xi \\ \eta \end{pmatrix} \right\rangle_{\R^{2n}} dt \\
     &=\int_0^1 \bigl\langle \xi , -2 (\p_\lambda P_\lambda) \dot \xi
      -(\p_\lambda \dot P_\lambda) \xi - P_\lambda (\p_\lambda P_\lambda) \xi
      -(\p_\lambda P_\lambda) P_\lambda \xi \bigr\rangle_{\R^{2n}} \\
     &\qquad -2 \bigl\langle \dot \xi + P_\lambda \xi ,
      (\p_\lambda P_\lambda) \xi \bigr\rangle_{\R^{2n}} \: dt \\
     &=\int_0^1 \bigl\langle \xi , - (\p_\lambda \dot P_\lambda) \xi
      -P_\lambda (\p_\lambda P_\lambda) \xi + (\p_\lambda P_\lambda) P_\lambda \xi 
      \bigr\rangle_{\R^{2n}} \: dt \\
     &=- \int_0^1 \bigl\langle \xi , 
      (\p_t \p_\lambda P_\lambda) \xi \bigr\rangle_{\R^{2n}} \: dt \\
     &=- \int_0^1 \bigl\langle \dot \xi , 
      (\p_\lambda P_\lambda) \xi \bigr\rangle_{\R^{2n}} \: dt
      -\int_0^1 \bigl\langle \xi , (\p_t \p_\lambda P_\lambda) \xi
      +(\p_\lambda P_\lambda) \dot \xi \bigr\rangle_{\R^{2n}} \: dt \\
     &=-\int_0^1 \p_t \bigl\langle \xi , 
      (\p_\lambda P_\lambda) \xi \bigr\rangle_{\R^{2n}} \: dt \\
     &=-\bigl\langle \xi(1) , 
      \bigl( \p_\lambda P_\lambda \bigr) (1) \: \xi(1) \bigr\rangle_{\R^{2n}}
      +\bigl\langle \xi(0) , 
      \bigl( \p_\lambda P_\lambda \bigr) (0) \: \xi(0) \bigr\rangle_{\R^{2n}}
      =0
    \end{split}
  \end{equation}
  where in the last equality we used the boundary conditions for $\xi$
  and the skewsymmetric family of matrices $P_\lambda$.
  We apply both identities to derive one more crucial result.
  Let $\xi \in Ker \: A_\lambda$ and use the symmetry of the projection operator $P_\lambda^\perp$
  to obtain
  \begin{equation} \label{eq:crucial-ident}
    \begin{split}
      &\bigl\langle 
      \xi , 
      P_\lambda^\perp \bigl( \p_\lambda A_\lambda \bigr) P_\lambda^\perp \: \xi 
      \big\rangle_{L^2} \\
     &=\bigl\langle \xi , 
      (\p_\lambda A_\lambda) \: \xi \big\rangle_{L^2}
      =\bigl\langle \xi , 
      (-\p_\lambda B_\lambda -\p_\lambda Q_\lambda) \: \xi \big\rangle_{L^2} \\
     &=\int_0^1 \left\langle
      \begin{pmatrix} \xi \\ \eta \end{pmatrix} ,
      \begin{pmatrix} -\p_\lambda B_\lambda -\p_\lambda Q_\lambda & 0 \\ 
      0 & 0 \end{pmatrix}
      \: \begin{pmatrix} \xi \\ \eta \end{pmatrix}
      \right\rangle_{\R^{2n}} dt \\
     &=\int_0^1 \left\langle
      \begin{pmatrix} \xi \\ \eta \end{pmatrix} ,
      \left( -U^{-\sigma} \bigl(\p_\lambda S_{\lambda,U}\bigr) U^\sigma
      +\begin{pmatrix} -\p_\lambda B_\lambda & \p_\lambda P_\lambda \\ 
      -\p_\lambda P_\lambda & 0 \end{pmatrix} \right)
      \begin{pmatrix} \xi \\ \eta \end{pmatrix} \right\rangle_{\R^{2n}}
      dt \\
     &=-\left\langle 
      \begin{pmatrix} \xi(0) \\ \eta(0) \end{pmatrix} ,
      \left( \int_0^1 \Psi_{\lambda,U}(t)^T \: 
      \bigl( \p_\lambda S_{\lambda,U} \bigr) (t)
      \Psi_{\lambda,U}(t) \: dt \right)
      \begin{pmatrix} \xi(0) \\ \eta(0) \end{pmatrix} 
      \right\rangle_{\R^{2n}} \\
     &=-\left\langle \Psi_{\lambda,U}(1)
      \begin{pmatrix} \xi(0) \\ \eta(0) \end{pmatrix} ,
      \hat{S}_1(\lambda) \: \Psi_{\lambda,U}(1)
      \begin{pmatrix} \xi(0) \\ \eta(0) \end{pmatrix}
      \right\rangle_{\R^{2n}} \\
     &=-\bigl\langle
      \zeta(0) , \hat{S}_1(\lambda) \: \zeta(0) \big\rangle_{\R^{2n}} \\
    \end{split}
  \end{equation}
  where we used (\ref{eq:identity2}) as well as orthogonality of 
  $U(t)$ in the fifth equality, 
  the identity (\ref{eq:identity1}) in the last but one equality and
  $\zeta(0)=(\xi (0),\eta (0))=(\xi_U (0),\eta_U (0)) \in Ker \:(\1-\Psi_{\lambda,U}(1))$ 
  in the last one.

    We are ready to prove (\ref{eq:RS=specflow}). 
  The first two equalities are by definition 
  of the spectral flow and the crossing operator
  \begin{equation*}
    \begin{split}
      \mu_{Spec}
     &( \{ A_\lambda \}_{\lambda \in [-T,T]} ) \\
     &=\sum_{\{ \lambda \mid Ker \: A_\lambda \not= \{ 0 \} \}}
      sign \: \Gamma(\{ A_\lambda \}_{\lambda \in [-T,T]} ,\lambda) 
      \mid_{ker \: A_\lambda} \\
     &=\sum_{\{ \lambda \mid Ker \: A_\lambda \not= \{ 0 \} \}}
      sign \: \langle \cdot , P_\lambda \bigl(\p_\lambda A_\lambda \bigr) P_\lambda 
      \cdot \rangle_{L^2} 
      \mid_{ker \: A_\lambda} \\
     &=\sum_{\{ \lambda \mid det (\1-\Psi_{\lambda,U}(1)) \not=0 \}}
      - sign \: \langle \cdot , \hat{S}_1(\lambda) \cdot \rangle_{\R^{2n}}
      \mid_{Ker \: (\1-\Psi_{\lambda,U} (1))} \\
     &=\sum_{\{ \lambda \mid Graph \Psi_{\lambda,U}(1) \cap \Delta \not=\{0\} \}}
      sign \: \Gamma \bigl( Graph \: \Psi_{\cdot,U} (1) ,\Delta, 
      \lambda \bigr)
     \mid_{\Delta \cap Graph \: \Psi_{\lambda,U}(1)} \\
     &=\mu_{RS} \bigl( Graph \: (\lambda \mapsto \Psi_{\lambda,U}(1)) ,\Delta \bigr) \\
     &=\mu_{RS} \bigl(\gamma_2,\Delta \bigr)
     \end{split}
  \end{equation*}
  where the third equality is (\ref{eq:crucial-ident}).
  Here it is important that the sums are over the same set of
  $\lambda's$. This follows from the equivalence of (S1) and (S4)
  $$
    \xi \in Ker \: A_\lambda \Longleftrightarrow
    \begin{pmatrix} \xi (0) \\ \dot \xi (0) +P_\lambda(0) \xi (0) \end{pmatrix}
    \in Ker \: \bigl( \1 - \Psi_{\lambda,U}(1) \bigr)
  $$
  Equality four has been derived in (\ref{eq:crossform=S}), 
  equality five is by definition of $\mu_{RS}$ and
  the last one is by definition of $\gamma_2$.
  This completes the proof of theorem \ref{thm:specflow}.
\end{proof}

\subsection{Proof of the index theorem \ref{thm:index}} \label{subsec:proofINDEX}
\begin{proof} Let $x \in Crit$ be a nondegenerate critical point,
i.e. $ker \: A_x =\{ 0 \}$,
and denote by $\mu_1 \le \mu_2 \le \ldots \le \mu_{Ind(x)} <0$
the negative eigenvalues of $A_x$ counted with multiplicities. 
Fix a real number $\hat{\mu} < \mu_1$.
It will be convenient in (\ref{eq:kappa}) to assume in addition $\hat{\mu} < -\pi$.
Let $\beta : [0,1] \to [\hat{\mu} ,0]$ be a smooth cut-off
function which equals $0$ near $0$ and $\hat{\mu}$ near $1$ and is strictly decreasing
elsewhere. For $i=1,\ldots,Ind(x)$ define $\lambda_i$ by $\beta (\mu_i)=\lambda_i$.
Let $Q$ be the path of symmetric matrices defined in (\ref{eq:A0})
and modify $\beta$, if necessary, such that its value at each $\lambda_i$ remains the same but
$\beta^\prime (\lambda_i) \notin Spec \: Q(1)$.
This technical condition ensures regularity of crossings.
\begin{figure}[ht]
  \centering\epsfig{figure=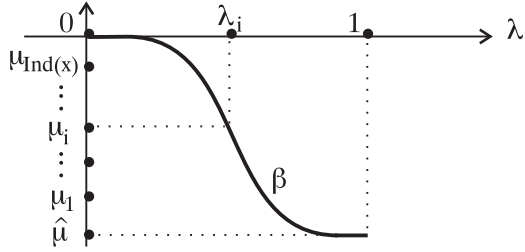}
  \caption{Cut-off function $\beta$ and negative spectrum of $A_x$} 
  \label{fig:cut-offf}
\end{figure}

  For $\lambda \in [0,1]$ consider two families of selfadjoint operators
in $L^2_\sigma$ with dense domain $W^{2,2}_\sigma$
\begin{equation} \label{eq:Ilambda-2}
  A_\lambda \xi
  = - \p_t \p_t \xi -B \xi - (1-\lambda) Q\xi -\beta(\lambda) \xi
\end{equation}
and
\begin{equation} \label{eq:Tilde-Ilambda-2}
  \tilde{A}_\lambda \xi
  = - \p_t \p_t \xi -B_\lambda \xi - \hat{\mu} \xi
\end{equation}
where 
$$
  B_\lambda = 2P_\lambda \p_t + (\p_t P_\lambda) + {P_\lambda}^2
  , \qquad P_\lambda = (1-\lambda) P
$$
and $B=B_0$. 
The reason to use two families instead of only one
is that regularity can be checked easily that way.
Both families fit into the framework
of section \ref{subsec:specflow-CZ}: Using the boundary conditions for $\xi$
and $\hat{\mu} <0$ a short calculation shows that the second family consists for
each $\lambda$ of a positive definite operator and so the injectivity
assumption (\ref{eq:injectivity-assumption}) is satisfied.
The family is also regular, because there are no crossings at all.
To check for the first family we observe that
$A_1=\tilde{A}_0$ is positive definite
and $A_0$ is precisely the model operator (\ref{eq:A0}) for $A_x$.
In view of $V \in \VMC_{reg}$ this confirms the injectivity 
assumption (\ref{eq:injectivity-assumption}).
Our choice of cut-off function now guarantees regularity of the family.
More precisely, assume $\lambda_i$ is a crossing and let
$\xi \in ker \: \p_\lambda A_{\lambda_i} \cap ker \:
A_{\lambda_i}$, then the first condition implies
\begin{equation} \label{eq:not-specQ}
  \begin{gathered}
    Q(t) \xi(t) = \beta^\prime (\lambda_i) \xi(t) \quad \forall t \in [0,1] \\
    \xi(0)=E_\sigma \xi(1) \; , \; \dot \xi(0) = E_\sigma \dot \xi(1)
  \end{gathered}
\end{equation}
Now $\beta^\prime (\lambda_i) \notin spec \: Q(1)$ implies $\xi(1)=0$.
Differentiating (\ref{eq:not-specQ}) with respect to $t$ leads
at $t=1$ to
$$
  ( \beta^\prime \smb{\lambda_i} - Q\smb{1} ) \dot \xi (1)
  = \dot Q(1) \xi(1) =0
$$
and in view of the regularity of the matrix $\beta^\prime (\lambda_i) - Q(1)$ we 
get $\dot \xi(1)=0$. Now we use the second condition which says that
$\xi$ is also in the kernel of the second order differential operator
$A_{\lambda_i}$. Since its boundary conditions are zero, it follows $\xi=0$.
This proves nondegeneracy of the crossing operator
$\Gamma( \{ A_\lambda \}_{\lambda \in [0,1]} , \lambda_i)$ at each crossing $\lambda_i$.

  We are in position to apply theorem \ref{thm:specflow} to both families.
Let us start with the first one and observe that its spectral
flow is given by $Ind(A_0)=Ind(x)$. We obtain
$$
  Ind(x) = \mu_{CZ} (\Psi_{1,U}) - \mu_{CZ} (\Psi_{0,U})
$$
where according to $(S4)$ in section \ref{subsec:specflow-CZ} for each $\lambda \in [0,1]$ 
the path $\Psi_{\lambda,U} : [0,1] \to Sp(2n)$ is determined by
\begin{equation*} 
  \p_t \Psi_{\lambda,U} = -J_0 S_{\lambda,U} \Psi_{\lambda,U}
  \; , \; \Psi_{\lambda,U} (0) = \1
\end{equation*}
and
\begin{equation*} 
  S_{\lambda,U} = U^\sigma
  \begin{pmatrix} (1-\lambda) Q +\beta(\lambda) \1 & P \\ 
  - P & \1 \end{pmatrix} U^{-\sigma}
  -J_0 U^\sigma \p_t (U^{-\sigma})
\end{equation*}
where $\sigma=\sigma(x) \in \{ 0,1 \}$.
Observe that $S_{0,U}$ is precisely the matrix $S_U$ in lemma (\ref{lem:IVP})
and so $\Psi_{0,U}$ equals $\gamma$ from (\ref{eq:linflow}).
In other words $\mu_{CZ} (\Psi_{0,U}) = \mu_{CZ} (z_x)$.

  Since $S_{1,U}$ contains the $t$-dependent matrices $P$ we make use
of the second family $\tilde{A}_\lambda$ in order to further simplify
the problem of calculating $\mu_{CZ}(\Psi_{1,U})$.
As before we get for each $\lambda \in [0,1]$ a path
$\tilde{\Psi}_{\lambda,U} : [0,1] \to Sp(2n)$ which is determined by
\begin{equation} \label{eq:PsiTilde}
  \p_t \Tilde{\Psi}_{\lambda,U} = -J_0 \Tilde{S}_{\lambda,U} \Tilde{\Psi}_{\lambda,U}
  \; , \; \Tilde{\Psi}_{\lambda,U} (0) = \1
\end{equation}
where
\begin{equation} \label{eq:STilde}
  \Tilde{S}_{\lambda,U} = U^\sigma
  \begin{pmatrix} \hat{\mu} \1 & (1-\lambda) P \\ 
  - (1-\lambda) P & \1 \end{pmatrix} U^{-\sigma}
  -J_0 U^\sigma \p_t (U^{-\sigma})
\end{equation}
The spectral flow of a family of positive definite
operators is zero and so theorem \ref{thm:specflow} gives
$$
  0= \mu_{CZ} (\Tilde{\Psi}_{1,U}) - \mu_{CZ} (\Tilde{\Psi}_{0,U}) 
  = \mu_{CZ} (\Tilde{\Psi}_{1,U}) - \mu_{CZ} (\Psi_{1,U})
$$
where we used $\tilde{S}_{0,U}=S_{1,U}$ in the last step.
Hence it remains to show that the Conley-Zehnder index of
the path $\Tilde{\Psi}_{1,U}$ equals $\sigma(x)$. In order to do so we
would like to treat the two cases $\sigma=0$ and $\sigma=1$ separately. \\
{\bf Case 1: \boldmath $\sigma=0$} We conclude

\begin{equation*}
  \begin{split}
   &\mu_{CZ} (\Tilde{\Psi}_{1,U}) \\
   &=\mu_{RS} (Graph \: \Tilde{\Psi}_{1,U} , \Delta ) \\
   &= {\textstyle \frac{1}{2}} sign \: \Gamma ( Graph \: \Tilde{\Psi}_{1,U} , \Delta , 0)
    = - {\textstyle \frac{1}{2}} sign \: \langle \cdot ,
    \Tilde{S}_{1,U}(0) \cdot \rangle \mid_{\R^{2n}} \\
   &= - {\textstyle \frac{1}{2}} sign \: \Big\langle \cdot , 
    \begin{pmatrix} \hat{\mu} \1 & 0 \\ 0 & \1 \end{pmatrix}
    \cdot \left. \Big\rangle \right|_{\R^{2n}} \\
   &= 0
  \end{split}
\end{equation*}
where we used (\ref{eq:RS=CZ}) in the first equality, (\ref{eq:crossform=S}) in the third
and the fact that there are no further crossings in the second equality.

\noindent
It remains to show that there are indeed no crossings $t>0$. This
is equivalent to $1$ not being in the spectrum of $\Tilde{\Psi}_{1,U}(t)$ 
for any $t\in (0,1]$.
Because the matrix
$$
  \Tilde{S}_{1,U} = \begin{pmatrix} \hat{\mu} \1 & 0 \\ 0 & \1 \end{pmatrix}
$$
is constant in $t$, we can integrate the corresponding differential equation (\ref{eq:PsiTilde})
for $\Tilde{\Psi}_{1,U}$ and obtain with $\kappa=\sqrt{-\hat{\mu}}$
$$
    \Tilde{\Psi}_{1,U}(t)
    =e^{-t J_0 \Tilde{S}_{1,U}}
    =exp \begin{pmatrix} 0&t \: \1 \\ -t \hat{\mu} \: \1 &0 \end{pmatrix}
    =\begin{pmatrix} \cosh t \kappa \: \1& 
    \kappa^{-1} \sinh t \kappa \: \1 \\ 
    \kappa \sinh t \kappa \: \1 &
    \cosh t \kappa \: \1 \end{pmatrix}
$$
We study the characteristic polynomial of $\Tilde{\Psi}_{1,U}(t)$ and obtain its eigenvalues
(of multiplicity $n$ each)
\[
  \rho_\pm(t)=\cosh t \kappa
  \pm \sinh t \kappa
\]
Finally the sum (difference) of the hyperbolic cosine and sine
is $1$ if and only if $t=0$, which shows that there are no further crossings.
Let us remark that, because $\Tilde{\Psi}_{1,U}(t)$ is symplectic, it follows
$\rho_+(t)=\rho_-(t)^{-1}$
and this is reflected in the key identity for hyperbolic
functions $\cosh^2 - \sinh^2 \equiv1$.

\noindent
{\bf Case 2: \boldmath $\sigma=1$} We obtain
\begin{equation*}
  \begin{split}
   &\mu_{CZ} (\Tilde{\Psi}_{1,U}) \\
   &=\mu_{RS} (Graph \: \Tilde{\Psi}_{1,U} , \Delta ) \\
   &= {\textstyle \frac{1}{2}} sign \: \Gamma ( Graph \: \Tilde{\Psi}_{1,U} , \Delta , 0)
    = - {\textstyle \frac{1}{2}} sign \: \langle \cdot , 
    \Tilde{S}_{1,U}(0) \cdot \rangle \mid_{\R^{2n}} \\
   &= - {\textstyle \frac{1}{2}} sign \: \Big\langle \cdot , 
    \begin{pmatrix} diag_n (\hat{\mu}-\pi, \hat{\mu} , \ldots , \hat{\mu} ) & 0_n \\ 
    0_n & diag_n (1-\pi, 1 , \ldots , 1) \end{pmatrix}
    \cdot \left. \Big\rangle \right|_{\R^{2n}} \\
   &=1
  \end{split}
\end{equation*}
where we used (\ref{eq:RS=CZ}) in the first equality, (\ref{eq:crossform=S}) in the third
and the fact that there are no further crossings in the second equality:
To see this let us rearrange coordinates $(x_1,\dots,x_n,y_1,\dots,y_n)$
of $\R^{2n}$ in the form $(x_1,y_1,x_2,\dots,x_n,y_2,\dots,y_n)$
such that $Sp(2n)$ gets identified with $Sp(2) \oplus Sp(2n-2)$ and
the path of symmetric matrices $\Tilde{S}_{1,U}(t)$ from (\ref{eq:STilde}) with
$$
  \begin{pmatrix} 
    \hat{\mu} \cos^2 \pi t + \sin^2 \pi t - \pi &
    (\hat{\mu} -1) \cos \pi t \sin \pi t \\
    (\hat{\mu} -1) \cos \pi t \sin \pi t &
    \hat{\mu} \sin^2 \pi t + \cos^2 \pi t - \pi
  \end{pmatrix}
  \oplus
  \begin{pmatrix} 
    \hat{\mu} \1_{n-1} &
    0_{n-1} \\ 
    0_{n-1} &
    \1_{n-1}
  \end{pmatrix}
$$
The symplectic path generated by the second term in the direct sum does not meet
the Maslov cycle for any $t\in (0,1]$ as was shown in case 1.
Let us denote the first term in the direct sum by $b(t)$.
It remains to investigate the path in $Sp(2)$ determined by
\begin{equation} \label{eq:gamma_2}
  \p_t \gamma_2(t) = -J_0 b(t) \gamma_2(t) , \qquad \gamma_2(0)=\1
\end{equation}
and show that there are no crossings with the Maslov cycle
$\CMC$ in $Sp(2)$ for any $t\in (0,1]$. A numerical plot of
$\gamma_2([0,1])$ confirming this
in case $\hat{\mu}=-\pi$ is shown in figure \ref{fig:mascyc} (b).
To give a proof we introduce the notation
$$
  u(t)=\begin{pmatrix} \cos \pi t & - \sin \pi t \\ \sin \pi t & \cos \pi t \end{pmatrix}
  , \qquad 
  s=\begin{pmatrix} \hat{\mu} & 0 \\ 0 & 1 \end{pmatrix}
$$
and observe that
$$
  b=usu^{-1}-J_0 u \p_t (u^{-1})
$$
Now $\gamma_2$ is a solution of (\ref{eq:gamma_2}) if and only if
$\psi=u^{-1} \gamma_2$ is a solution of
$$
  \p_t \psi =-J_0 s \psi , \qquad \psi(0) = \1
$$
where we used the identity $\pi J_0= (\p_t u) u^{-1}$.
The solution $\psi (t)$ is known from case $1$ above.
We are done once we have shown that the following function
is strictly positive for $t \in (0,1]$
\begin{equation*}
  \begin{split}
    f(t)
   &=\det \bigl( \1 - \gamma_2(t) \bigr)
    =\det \bigl( u^{-1}(t) - \psi(t) \bigr) \\
   &=2-2 \cos \pi t \cosh \kappa t
    +(\kappa - \kappa^{-1}) \sin \pi t \sinh \kappa t
  \end{split}
\end{equation*}
Note that $f(0)=0$ and $f(1)=2(1+\cosh \kappa) >2$.
The function is clearly positive for $t \in [1/2,1]$ if $\kappa \ge 1$.
This is true since
\begin{equation} \label{eq:kappa}
  \kappa = \sqrt{-\hat{\mu}} > \sqrt{\pi}
\end{equation}
We calculate the first derivative
$$
  \dot f (t)
  =(2\pi -1 +\kappa^2 ) \sin \pi t \cosh \kappa t
  +\bigl( (\pi -2) \kappa - {\textstyle \frac{\pi}{\kappa}} \bigr)
  \cos \pi t \sinh \kappa t
$$
which is positive on $(0,1/2)$ since both coefficients are, due
to (\ref{eq:kappa}).
Moreover, we get $\dot f (0)=0$ and $\ddot f (0)=\pi (2\pi -1 +\kappa^2 )
+\kappa ( (\pi -2) \kappa - \frac{\pi}{\kappa} ) >0$
so that $0$ is a local minimum and it follows that $f$ is strictly increasing on $(0,1/2)$.
This finishes the proof of case $2$ and of the index theorem \ref{thm:index}.
\end{proof}

\section{Transversality} \label{sec:transversality}

\subsection{Thom-Smale transversality} \label{subsec:thom-smale}

  We recall the basics of Thom-Smale transversality theory.
Let $\AMC, \BMC$ be smooth Banach manifolds which admit a countable
atlas each and are modeled on separable Banach spaces.
Let $\EMC \to \AMC \times \BMC$ denote a smooth Banach space bundle,
$\EMC_{(a,b)}$ the fibre over $(a,b)$
and $\FMC$ a section of $\EMC$ of class $C^\ell$, $\ell \ge 1$.
Define $\FMC_b(a)=\FMC_a(b)=\FMC(a,b)$ for
$a\in \AMC$, $b\in \BMC$, and let $d$ denote the differential
of a section followed by projection onto the fibre.

\begin{thm} \label{thm:Thom-Smale}
  Let $\FMC$ be a $C^\ell$-section of $\EMC$ as above
  and assume

  \vspace{2mm}
  \noindent
  \begin{minipage}{.8cm}
    $(F)$
  \end{minipage}
  \;
  \begin{minipage}{11.2cm} \it
    $d\FMC_b (a) : T_a \AMC \to \EMC_{(a,b)}$ is Fredholm for all $a \in \FMC_b^{-1} (0)$
    and $\ell \ge \max \{ 1,1+Ind \: d\FMC_b (a) \}$
  \end{minipage}

  \vspace{2mm}
  \noindent
  \begin{minipage}{.8cm}
    $(S)$
  \end{minipage}
  \;
  \begin{minipage}{11.2cm} \it
    $d\FMC (a,b) : T_a \AMC \times T_b \BMC \to \EMC_{(a,b)}$ 
    is surjective for all $(a,b) \in \FMC^{-1} (0)$
  \end{minipage}
  \vspace{2mm}
  \noindent

  \noindent
  In this case the subset
  $$
    \{ b\in \BMC \mid d\FMC_b (a) \; \mbox{is surjective for all} \; a\in \FMC_b^{-1} (0) \}
    \subset \BMC
  $$
  is residual and therefore dense.
\end{thm}

\begin{defn}
  Let $W$ be a closed subspace of a Banach space $X$.
  A closed subspace $V \subset X$ is called a \emph{topological complement of} $X$,
  if $W \oplus V =X$ and $W \cap V = \{ 0 \}$.
\end{defn}

One crucial ingredient in the proof of theorem \ref{thm:Thom-Smale} 
is the following well known proposition.
For the sake of completeness we give a proof following \cite{Sa96}.

\begin{prop} \label{prop:folk}
  Let $X,Y,Z$ be Banach spaces. Assume $D:X\to Y$ is Fredholm
  and $L:Z\to Y$ is a bounded linear operator, then

  \noindent
  $i)$ the range of $D\oplus L : X\oplus Z \to Y$
    is closed with a finite-dimensional complement

  \noindent
  $ii)$ if $D\oplus L$ is surjective, then $ker \: (D\oplus L)$ admits a topological
    complement.
    Moreover, the projection on the second factor
    $
      \Pi : ker \: (D\oplus L) \to Z 
    $
    is Fredholm with
    $$
      ker \: \Pi \simeq ker \: D , \qquad coker \: \Pi \simeq coker \: D
    $$
    where $coker \: D := Y / ran \; D$ and similarly for $coker \: \Pi$.
\end{prop}

\begin{proof}
  {\bf ad \boldmath $i)$}
  Since $ran \; D$ is closed the quotient $Y/ran \; D$ inherits the structure
  of a Banach space from $Y$ and the projection $pr : Y \to Y/ran \; D$ is continous.
  By assumption $Y/ran \; D$ is finite dimensional and so is its subspace
  $pr ( ran \; L)$. Hence the latter space is closed and so is its preimage
  $pr^{-1}(pr ( ran \; L))$
  under the continous map $pr$, which is
  $ran \: D + ran \; L$. 
  Finally, the complement of $ran \; (D\oplus L)$ in $Y$ is contained
  in the finite dimensional space $Y/ran \; D$.
  \\
  {\bf ad \boldmath $ii)$}
  It is well known, cf. \cite{Br83} section II.4, that any finite dimensional
  subspace of a Banach space admits a topological complement.
  Let $X_1$ be a topological complement of $ker \: D$.
  Moreover, since $ran \; D$ is closed with finite dimensional complement, we
  can write $Y=ran \; D \oplus coker \; D$.
  Surjectivity of $D\oplus L$ implies that $coker \; D \subset ran \; L$
  and so we may choose a basis $\{ Lz_1 , \ldots , Lz_N \}$ of $coker \; D$,
  where $\{ z_1 , \ldots , z_N \}$ is a set of linearly independent elements of $Z$.
  \\
  Now define the linear operator 
  $T : ran \: D \oplus coker \: D \to ker \: D \oplus X_1 \oplus Z$
  $$
    T(y^\prime,y^{\prime\prime})= (0,x^\prime,\sum_{\nu =1}^{N} \lambda_\nu z_\nu )
  $$
  where $x^\prime$ is determined uniquely by $y^\prime=Dx^\prime$ and
  the coefficients $\lambda_\nu$ by $y^{\prime\prime} =\sum_{\nu =1}^{N} \lambda_\nu L z_\nu$.
  It is not hard to check that a topological complement of $ker \: (D\oplus L)$
  is given by $ran \; T$. Actually, since existence of a right inverse of a bounded linear
  operator is equivalent to surjectivity and existence of a topological complement
  of the kernel, we observe that $T$ is indeed a right inverse of $D\oplus L$.
  \\
  Now assume $(x,z) \in ker \: (D \oplus L)$, then
  $$
    (x,z) \in ker \: \Pi
    \Leftrightarrow Dx=0 \; \text{and} \; z=0
    \Leftrightarrow (x,0) \in ker \: D \oplus \{ 0 \}
  $$
  which proves that $ker \: \Pi =ker \: D \oplus 0$.
  Define
  \[
    L^{-1} (ran \: D) := 
    \{ z\in Z \mid Lz =Dx \quad \mbox{for some} \quad x \in X \}
  \]
  then $ran \: \Pi =L^{-1} (ran \:D)$ and this set is closed:
  $ran \: D$ is closed and so is its preimage under the continuous map $L$.
  Finally we obtain
  \[
    coker \: \Pi
    = \frac{Z}{ran \: \Pi}
    = \frac{Z}{L^{-1}(ran \: D)}
    \simeq \frac{ran \: L}{ran \: D \cap ran \: L}
    = \frac{Y}{ran \: D}
    = coker \: D
  \]
  where the first isomorphism is induced by $L$, namely
  \begin{equation*}
    \begin{split}
      \frac{Z}{L^{-1}(ran \: D)}
     &\to \frac{ran \: L}{ran \: D \cap ran \: L} \\
      [ z+L^{-1}(ran \: D) ] 
     &\mapsto [ Lz+ ran \: D \cap ran \: L ]
    \end{split}
  \end{equation*}
  and the subsequent identity is due to
  $Y=ran \: D + ran \: L$ and given by
  \begin{equation*}
    \begin{split}
      \frac{ran \: D + ran \: L}{ran \: D}
     &\to \frac{ran \: L}{ran \: D \cap ran \: L} \\
      [ Lz+ ran \: D \cap ran \: L ] 
     &\mapsto [ Lz+ ran \: D ]
    \end{split}
  \end{equation*}
\end{proof}

  Since  $d\FMC_a (b)$ is bounded ($\ell \ge1$) and
$d\FMC (a,b) =  d\FMC_b (a) \oplus  d\FMC_a (b)$,
part $i)$ of proposition \ref{prop:folk} tells that $ran \: d\FMC (a,b)$ is closed. Therefore
in order to verify $(S)$ it is sufficient to prove its density
for all $(a,b) \in \FMC^{-1}(0)$. Now it is a consequence
\cite[cor. I.8]{Br83} of the Hahn-Banach theorem that the latter
is equivalent to its annihilator being trivial

\vspace{3mm}
\noindent

\noindent
\begin{minipage}{.8cm}
  $(A)$
\end{minipage}
\;
\begin{minipage}{11.2cm} \it
  $\{ v^* \in \EMC_{(a,b)}^* \mid v^* (v)=0 \; \forall v \in ran \: d\FMC (a,b) \} =
   \{ 0 \} \; \; \forall (a,b) \in \FMC^{-1} (0)$
\end{minipage}

\vspace{3mm}

\noindent
In many applications it is more convenient to check conditions
$(F)$ and $(A)$  instead of $(F)$ and $(S)$.

  We sketch the proof of theorem \ref{thm:Thom-Smale}: Properties $(F)$ and $(S)$ together
with part $ii)$ of proposition \ref{prop:folk} tell that for any $(a,b) \in \FMC^{-1} (0)$
its linearization $d\FMC (a,b)$ admits a topological complement
and therefore a right inverse. This means, by definition, that $0$ is
a regular value of $\FMC$ and so, by the implicit function theorem,
the so called \emph{universal moduli space}
$$
  X=\FMC^{-1}(0)
$$
is a $C^\ell$-Banach manifold. It is locally at $(a,b)$ modeled on the separable
Banach space $ker \: d\FMC (a,b)$ and admits a countable atlas.
Define the projection onto the second factor
$\pi : \AMC \times \BMC \supset X \to \BMC$ and observe that, by
$ii)$ of the Lemma,
$$
  d\pi (a,b) : ker \: \left( d\FMC_b (a) \oplus  d\FMC_a (b) \right)
  \to T_b \BMC , \qquad (A,B) \mapsto B
$$
is Fredholm for any $(a,b)\in \FMC^{-1}(0)$ with
$$
  Ind \: d\pi (a,b) = Ind \: d\FMC_b (a)
$$ 
This means, again by definition, that $\pi$ is a Fredholm map of class $C^\ell$ between
separable Banach manifolds.
Finally we may apply the Sard-Smale theorem \cite{Sm73}
to the countably many coordinate representations of $\pi$
and obtain that the set of regular values of $\pi$ is
a residual subset of $\BMC$ in case $\ell \ge \max
\{ 1,1+Ind \: d\FMC_b(a) \}$.
Now the theorem follows from the following well known result of Thom-Smale
transversality theory

\begin{lem} \label{lem:wonder-TS}
  Under the assumptions of theorem \ref{thm:Thom-Smale} it follows
  $$
    \{ \mbox{regular values of} \; \pi \}
    =\{ b\in \BMC \mid d\FMC_b(a) \; \mbox{is surjective} \; \forall a\in \FMC_b^{-1}(0) \}
  $$
\end{lem}

\begin{proof}
  We have to show that for every $b \in \BMC$ it holds
  $$
    d\pi (a,b) \; \text{surjective} \; \forall a \in \FMC_b^{-1} (0)
    \quad \Longleftrightarrow \quad 
    d \FMC_b (a) \; \text{surjective} \; \forall a \in \FMC_b^{-1} (0)
  $$
  Assume the right hand side was true and choose any $\hat{b} \in T_b\BMC$.
  By surjectivity of $d \FMC_b(a)$ we find $\hat{a} \in T_a\AMC$ such that
  $d \FMC_b(a) \hat{a} = d \FMC_a(b) \hat{b}$. Then the pair 
  $(-\hat{a},\hat{b}) \in T_a\AMC \times T_b\BMC$
  is indeed element of $T_{(a,b)}X$ and is mapped to $\hat{b}$
  under $d\pi (a,b)$. \\
  To prove the other direction assume the left hand side was true,
  choose any $v \in \EMC_{(a,b)}$ and apply the surjectivity assumption $(S)$
  of theorem \ref{thm:Thom-Smale} to conclude the existence of a pair
  $(\hat{a_0},\hat{b}) \in T_a\AMC \times T_b \BMC$ being mapped to $v$
  by $d\FMC(a,b)$; in other words
  $$
    d \FMC_b(a) \hat{a_0} + d \FMC_a(b) \hat{b} = v
  $$
  Given $\hat{b}$, surjectivity of $d \pi (a,b)$ implies existence
  of $\hat{a_1}$ such that $(\hat{a_1},\hat{b}) \in T_{(a,b)}X$, i.e.
  $d \FMC_b(a) \hat{a_1} + d \FMC_a(b) \hat{b} = 0$.
  Now define $\hat{a}=\hat{a_0}-\hat{a_1}$ and observe that 
  $d \FMC_b(a) (\hat{a_0} - \hat{a_1})=v$.  
\end{proof}

\subsection{Transversality in the $C^k$-category} 
  \label{subsec:transversality-Ck}

  In this subsection we fix $k\ge2$ and apply Thom-Smale transversality
theory to the following situation
$$
  \AMC=W^{2,2}(S^1,M) , \qquad \BMC=C^k(S^1\times M,\R)
$$
The $C^{k-1}$-section  $\FMC : \AMC \times \BMC \to \EMC$
is defined to be
\begin{equation} \label{eq:F}
  \FMC (x,V) := grad \: \SMC_V (x) = -\n_t \dot x - \n V_t(x)
\end{equation}
The fibre of $\EMC$ at $(x,V)$ is $\EMC_{(x,V)}=L^2(S^1,x^*TM)$ and
the linearization of $\FMC$ at a zero turns out to be
\begin{equation*}
  \begin{gathered}
    d\FMC(x,V) : W^{2,2}(S^1,x^*TM) \times C^k(S^1\times M,\R) 
    \to L^2(S^1,x^*TM) \\
    (\xi,\hat{V}) 
    \mapsto d\FMC_V(x) \xi + d\FMC_x(V) \hat{V}
  \end{gathered}
\end{equation*}
where
\begin{equation*}
  \begin{split}
    d\FMC_V(x) \xi
   &=-\n_t \n_t \xi -R(\xi,\dot x)\dot x -\n_\xi \n V_t(x) 
    = A_x \\
    d\FMC_x(V) \hat{V}
   &=-\n \hat{V}_t(x) \\
  \end{split}
\end{equation*}

\begin{rmk} \label{rmk:separabel}
  1) Consider the metric space $W^{2,2}(S^1,\R^N)$. It is separabel
  because the set of Fourier-series with coefficients in $\Q^N$ is countable
  and dense in $C^\infty(S^1,\R^N)$, which in turn is dense
  in $W^{2,2}(S^1,\R^N)$.

  2) Recall that any metric space $X$ is paracompact, which by definition means that
  any open cover admits a locally finite refinement.
  If $X$ is in addition separable, we can conclude that any open
  cover admits a countable subcover. The following argument is taken from
  \cite{Sa96}: if the cover $\{ U_\alpha \}_\alpha$ is locally finite and
  $\{ x_i \}_i$ is a dense sequence then the set of pairs
  $(\alpha,i)$ with $x_i \in U_\alpha$ is countable. Since every $U_\alpha$
  contains some point $x_i$ the map $(\alpha,i) \mapsto \alpha$ is surjective.

  3) For sufficiently large $N\in \N$ we can find an embedding
  $\iota : M \hookrightarrow \R^N$ and define $W^{2,2}(S^1,M)$ to be the set of those elements
  of $W^{2,2}(S^1,\R^N)$ whose image lies in $\iota(M)$.
  That way $W^{2,2}(S^1,M)$ inherits the metric from the ambient space
  and is therefore itself paracompact and separabel.

  4) One can give $W^{2,2}(S^1,M)$ the structure of a Banach manifold
  modeled on the separabel Banach space $W^{2,2}(S^1,\R^n)$, where $n=dim \: M$. 
  An atlas can be constructed where the charts
  are labeled by smooth functions from $S^1$ to $M$. Since $W^{2,2}(S^1,M)$
  is paracompact and separabel by 3), we can apply 2) to get a countable
  subatlas.
\end{rmk}

\begin{thm} \label{thm:Ck-transversality}
  {\rm ({\bf Transversality in $C^k$})}
  For every integer $k\ge 2$ the following are true:

  \noindent
  $i)$ The subset
  $$
    \VMC^k_{reg} =\{ V \in C^k(S^1 \times M,\R) \mid
    d\FMC_V(x) \; \mbox{is surjective} \; \forall x \in \FMC_V^{-1}(0) \}
  $$
  of the Banach space $( C^k(S^1 \times M,\R) , \| \cdot \|_{C^k} )$
  is residual and therefore dense.

  \noindent
  $ii)$ Fix $a\in \R$, then the set
  $$
    \VMC^{k,a}_{reg} =\{ V \in C^k(S^1 \times M,\R) \mid
    d\FMC_V(x) \; \mbox{is surjective} \; \forall x \in \FMC_V^{-1}(0),
    \SMC_V (x) < a \}
  $$
  is open and dense in $(C^k(S^1 \times M,\R) , \| \cdot \|_{C^k})$.
\end{thm}

\begin{proof}
  {\bf ad \boldmath $i)$} We prove that assumption $(F)$ holds. Let
  $x \in \FMC_V^{-1}(0)$. Recall that
  $d\FMC_V(x)$ is the perturbed Jacobi-operator $A_x$ and
  $$
    dim \: ker \: d\FMC_V(x) = dim \: ker \: A_x
    = Nullity (x) < \infty
  $$
  by the Morse index theorem. Moreover, since $d\FMC_V(x)$ is selfadjoint
  $$
    coker \: d\FMC_V(x) \simeq ker \: d\FMC_V(x) 
  $$
  This shows that $d\FMC_V(x)$ is Fredholm of index $0$ for any
  $x \in \FMC_V^{-1}(0)$. 
  Together with the assumption $k\ge2$ it follows that
  the condition in Theorem \ref{thm:Thom-Smale}
  on the differentiability of $\FMC$ is satisfied: $k-1 \ge 1$.

    It remains to verify assumption $(A)$, i.e. we have to show that
  \[
    \langle \eta , d\FMC_V (x) \xi \rangle =0 \quad \forall \xi \in W^{2,2}(x^*TM)
  \]
  and
  \[
    \langle \eta , d\FMC_x (V) \hat{V} \rangle =0 \quad \forall \hat{V} \in 
    C^k(S^1 \times M,\R)
  \]
  together imply $\eta =0$. The first condition says that $\eta \in ker\: d\FMC_V (x)$.
  So it satisfies a linear second order ODE with coefficients of class $C^{k-2}$
  and therefore $\eta \in C^k(x^*TM)$. Now assume by contradiction that there
  is $t_0 \in S^1$ such that $\eta(t_0) \not= 0$. In five steps we are going to construct
  $\hat{V}_t \in C^\infty$  such that
  \[
    \langle \eta , \n \hat{V}_t (x) \rangle_{L^2} > 0 .
  \]
  
    As our construction will be local, we may choose geodesic
  normal coordinates $\xi=(\xi^1,\ldots,\xi^n)$ around $x_0=x(t_0)$.
  Let $\iota$ denote the injectivity radius of $(M,g)$. The piece of the 
  loop $x(t)$ which lies inside the
  coordinate patch determines the curve $\xi(t)$ via
  \[
    x(t) =exp_{x_0} \: \xi(t)
  \]
  so that $\xi(t_0)=0$.
  
  {\sc Step 1}
    Because $x(t)$ is continuous, we may choose a constant $\delta_1 >0$
    sufficiently small such that
    \[
      \mid \xi (t) \mid \le \iota / 2 , \qquad \forall
      t \in [ t_0-\delta_1 , t_0+\delta_1 ] .
    \]
  
  {\sc Step 2}
    Because $\eta$ is continuous and $\eta (t_0) \not= 0$, we may choose a
    constant $\delta_2 >0$ sufficiently small such that
    \[
      \langle \eta (t) , \eta (t_0) \rangle >0
      , \qquad \forall t \in [ t_0-\delta_2 , t_0+\delta_2 ] .
    \]
  
  {\sc Step 3}
    Set $\delta=\min \{ \delta_1 , \delta_2 \}$ and choose a cut-off function
    $\gamma \in C^\infty (\R ,[0,1])$ such that
    \[
      \gamma(t)=
      \begin{cases}
        1& , t \in [ t_0-\frac{\delta}{2} , t_0+\frac{\delta}{2} ] \\
        0& , t \notin [ t_0-\delta , t_0+\delta ] .
      \end{cases}
    \]
  
  {\sc Step 4}
    Choose a cut-off function $\beta \in C^\infty (\R ,[0,1])$ such that
    \[
      \beta(\mid \xi \mid^2)=
      \begin{cases}
        1& , \mid \xi \mid^2 \le \iota^2/2 \\
        0& , \mid \xi \mid^2 \ge \iota^2 .
      \end{cases}
    \]

  {\sc Step 5}
    We are ready to define $\hat{V}_t$
    \[
      \hat{V}_t (x) =
      \begin{cases}
        \gamma(t) \: \beta(\mid \xi \mid^2) \: 
        \langle \eta(t_0) , \xi \rangle
       & , x= exp_{x_0} \xi \; \mbox{and} \mid \xi \mid^2 < \iota^2 \\
        0
       & , \mbox{else} .
      \end{cases}
    \]

  \noindent
  Putting everything together we get
  \begin{equation*}
    \begin{split}
      \langle \eta , \n \hat{V}_t (x) \rangle_{L^2}
     &=\int_0^1  \langle \eta (t) , \n \hat{V}_t (x(t)) \rangle \: dt \\
     &=\int_0^1 d \hat{V}_t (x(t)) \comp \eta (t) \: dt \\
     &=\int_{ \{ t : |\xi(t)| < \iota \} }
      \frac{\p \hat{V}_t}{\p \xi^j} \mid_{exp_{x_0} \xi (t)} \: \eta^j(t) \: dt \\
     &=\int_{t_0-\delta}^{t_0+\delta} \Bigl(
      2 \gamma(t) \: \beta^\prime (\mid \xi \mid^2) \: 
      \langle \xi(t) , \eta (t) \rangle \:
      \langle \eta (t_0) , \xi (t) \rangle \\
     &\qquad \qquad \quad + \gamma(t) \: \beta (\mid \xi \mid^2) \:
      \langle \eta (t_0) , \eta (t) \rangle \Bigr) dt \\
     &=\int_{t_0-\delta}^{t_0+\delta} 
      \gamma(t) \: \langle \eta (t_0) , \eta (t) \rangle \: dt >0 .
    \end{split}
  \end{equation*}
  The third equality follows from the definition of $\hat{V}_t$ (Step 5),
  and the fourth one from Step 3 ($supp \: \gamma$) as well as a
  straight forward calculation. In the fifth equality we
  used that for $t \in [t_0-\delta , t_0+\delta]$ Step 1 implies
  $\mid \xi (t) \mid^2 \le \iota^2/4$
  and therefore, by Step 4, $\beta^\prime \equiv 0$ and $\beta \equiv 1$.
  Step 2 gives the final strict inequality.

    {\bf ad \boldmath $ii)$} \emph{dense:} 
  Since 
  $$
    \VMC^k_{reg} \subset \VMC^{k,a}_{reg}
  $$
  and $\VMC^k_{reg}$
  is a residual subset of $C^k(S^1 \times M,\R)$ by $i)$, the set $\VMC^{k,a}_{reg}$ is too.
  By Baire's category theorem it is therefore dense.

    \emph{open:} We fix some regular $V$ and construct an open neighbourhood
  $W_V$ of $V$ in $C^k(S^1 \times M,\R)$ such that $W_V \subset \VMC^{k,a}_{reg}$.
  Recall that
  $$
    \pi : X=\FMC^{-1}(0) \to C^k(S^1 \times M,\R)
  $$
  is a Fredholm map of class $C^{k-1}$. Since
  $$
    X^a = \{ (x,v) \in X \mid \SMC_V(x) < a \}
  $$
  is open in $X$, it follows that the restriction
  $$
    \pi_a : X^a \to C^k(S^1 \times M,\R)
  $$
  is Fredholm of class $C^{k-1}$ too. Now rewrite $\pi_a^{-1} (V)$ as
  $$
     \{ x\in W^{2,2}(S^1,M) \mid d\FMC_V (x) \; \mbox{surjective} , \;
      \SMC_V(x) < a \; \forall x\in \FMC_V^{-1}(0) \} \times \{ V \}
  $$
  and observe that this set is compact (cf. our remarks in the introduction
  on Cieliebak's uniform $C^0$-bound); in other words $\pi_a$ is a proper map.
  As will be discussed later surjectivity is an open condition and so
  there exists an open neighbourhood $U$ of $\pi_a^{-1}(V)$
  in $X^a$ such that $d\FMC_{V^\prime} (x^\prime)$ is surjective
  for all $(x^\prime , V^\prime ) \in U$. 
  Now use continouity of $\pi_a$ and compactness of $\pi_a^{-1} (V)$
  to conclude the existence of an open neighbourhood 
  $W_V$ of $V$ in $C^k(S^1 \times M,\R)$ such that
  $\pi_a^{-1} (W_V) \subset U$. It follows $W_V \subset \VMC^{k,a}_{reg}$.

    It remains to show that surjectivity of $d\FMC_V(x)$ is an open condition
  (in $X$): Note that for $(x^\prime,V^\prime)$ near $(x,V)$ the
  operators $d\FMC_V(x)$ and $d\FMC_{V^\prime} (x^\prime)$
  (strictly speaking their representatives with respect to
  a trivialization of $\EMC \to W^{2,2} \times C^k$ at $(x,V)$) 
  differ by a bounded operator, whose norm can be made arbitrarily small
  by choosing $(x^\prime,V^\prime)$ sufficiently close to $(x,V)$.
  In view of the subsequent lemma and the selfadjointness of the
  linear operators we get
  \begin{equation*}
    \begin{split}
      dim \: coker \: d\FMC_{V^\prime} (x^\prime)
     &=dim \: ker \: d\FMC_{V^\prime} (x^\prime) \\
     &\le dim \: ker \: d\FMC_{V} (x) \\
     &=dim \: coker \: d\FMC_{V} (x) =0 \\
    \end{split}
  \end{equation*}
\end{proof}

\begin{lem} \label{lem:kernel-decrease}
  Let $X,Y$ be Banach spaces and $D:X\to Y$ be a Fredholm operator.
  Then there exists an $\epsilon >0$ such that for any linear operator
  $L:X\to Y$ with $\| L \| < \eps$
  \[
    dim \: ker \: (D+L) \le dim \: ker \: D
  \]
\end{lem}
\
\begin{proof}
  Let $X_1$ be a topological complement of $ker \: D$.
  It suffices to show
  $$
    ker \: (D+L) \cap X_1 =\{ 0 \}
  $$
  The restriction $\tilde{D} : X_1 \to ran \: D$ is a bounded bijection
  and therefore admits a bounded inverse $\tilde{D}^{-1}$ by the open
  mapping theorem. Let $x\in X_1 \cap ker \: (D+L)$, i.e.
  $x=-\tilde{D}Lx$, then
  $$
    \| x \|_X=\| \tilde{D}^{-1} L x \|_X
    \le \| \tilde{D}^{-1} \| \cdot \| L \| \cdot \| x \|_X
    \le \eps \| \tilde{D}^{-1} \| \cdot \| x \|_X .
  $$
  Choose $0<\eps<\| \tilde{D}^{-1} \|$, then it follows $x=0$.
\end{proof}

\subsection{Transversality in the $C^\infty$-category} 
  \label{subsec:transversality-Cinfty} \unboldmath

\begin{proof} ({\sc of Theorem \ref{thm:transversality}})
{\bf ad \boldmath $i)$} $\VMC^a_{reg} \subset (C^\infty,d)$ \emph{dense:}
Given any $V \in C^\infty(M\times S^1,\R)$
we have to construct a sequence $V_k^\prime \in \VMC^a_{reg}$
such that 
\[
  \forall \eps >0 \quad \exists k_0 \in \N \quad \forall k>k_0
  \quad : \quad d(V,V_k^\prime) < \eps .
\]
The idea will be to approximate $V$ by regular $V_k$'s in the $C^k$-topology
and then approximate $V_k$ by smooth regular elements 
$V_k^\prime$ in the $C^k$-topology.
Finally we make use of the observation that in order to control the metric $d$
we essentially have to control only finitely many $C^k$-norms in its series,
because the strong weights $1/2^{k}$ take care of all the other ones.

  {\sc Step 1} Because $\VMC^{k,a}_{reg}$ is dense in
$(C^k\smb{M\times S^1,\R},\| \cdot \|_{C^k})$ for any integer $k\ge 2$,
we can find $V_k \in \VMC^{k,a}_{reg}$ with
$\| V-V_k \|_{C^k} < 1/(2k)$. For $k=0,1$ let us define $V_0=V_1=V$.

  {\sc Step 2} Because $\VMC^{k,a}_{reg}$ is open in
$(C^k\smb{M\times S^1,\R},\| \cdot \|_{C^k})$ for any integer $k\ge 2$,
we can choose $0<\eps_k <1/(2k)$ sufficiently small such that
$B_{\eps_k}(V_k)$, the open $\eps_k$-ball around $V_k$, is contained in
$\VMC^{k,a}_{reg}$.

  {\sc Step 3} Because $M\times S^1$ is compact, $C^\infty(M\times S^1,\R)$
is dense in $(C^k\smb{M\times S^1,\R},\| \cdot \|_{C^k})$ for $k \in \N_0$,
cf. \cite[theorem 2.6]{Hi76}. Hence we can find $V_k^\prime \in
C^\infty \cap B_{\eps_k}(V_k)$ for $k\ge 2$. For $k=0,1$ we define
$V_0^\prime=V_1^\prime=V$.

  {\sc Step 4} Now pick $\eps>0$ and choose $\nu_0 \in \N$ sufficiently large such that
$f(\nu_0)=\sum_{\nu =\nu_0 +1}^\infty 2^{-\nu} < \eps /2$.
Choose $k_0 > \max \{ \nu_0 , 4/\eps \}$ and observe that by Steps 1,2 and 3 
for $k\ge 2$
\[
  \| V-V_k^\prime \|_{C^k} \le \| V-V_k \|_{C^k}
  + \| V_k - V_k^\prime \|_{C^k}
  \le \frac{1}{2k} + \eps_k \le \frac{1}{k} .
\]
Note that this implies $\| V-V_k^\prime \|_{C^\nu} \le \| V-V_k^\prime \|_{C^k}
\le \frac{1}{k}$ for any $0 \le \nu \le k$.
We get for any $k>k_0$
\begin{equation*}
  \begin{split}
    d(V,V_k^\prime)
   &=\sum_{\nu =0}^{\nu_0} \frac{1}{2^\nu} \: 
    \frac{\| V-V_k^\prime \|_{C^\nu}}{1+\| V-V_k^\prime \|_{C^\nu}}
    +\sum_{\nu =\nu_0 +1}^\infty \frac{1}{2^\nu} \: 
    \frac{\| V-V_k^\prime \|_{C^\nu}}{1+\| V-V_k^\prime \|_{C^\nu}} \\
   &\le \sum_{\nu =0}^{\nu_0} \frac{1}{2^\nu} \: \frac{1}{k} \:
    \frac{1}{1+\| V-V_k^\prime \|_{C^\nu}}
    +\sum_{\nu =\nu_0 +1}^\infty \frac{1}{2^\nu} \\
   &\le \frac{2}{k} + \frac{\eps}{2} < \eps .
  \end{split}
\end{equation*}

  $\VMC^a_{reg} \subset (C^\infty,d)$ \emph{open:}
Pick $V \in \VMC^a_{reg}$ and set $k=2$.
Exploiting openess of $\VMC^{2,a}_{reg}$ in 
$(C^2\smb{M\times S^1,\R} ,\| \cdot \|_{C^2})$ we are able to choose
a constant $\eps_0 >0$ such that for any $V^{\prime\prime}$ of class $C^2$
with $\| V-V^{\prime\prime} \|_{C^2} < \eps_0$ it follows
$V^{\prime\prime} \in \VMC^{2,a}_{reg}$.
Now define
\[
  \eps = \frac{1}{4} \: \frac{\eps_0}{1+\eps_0} .
\]
Let $V^\prime$ of class $C^\infty$ be such that $d(V,V^\prime) <\eps$.
Therefore each term in the series on the left hand side
has to be strictly smaller then $\eps$,
in particular the second one
\[
  \frac{1}{2^2} \: \frac{\| V-V^\prime \|_{C^2}}{1+ \| V-V^\prime \|_{C^2}}
  < \eps =\frac{1}{4} \: \frac{\eps_0}{1+\eps_0} .
\]
But this is equivalent to
\[
  \| V-V^\prime \|_{C^2} < \eps_0
\]
and therefore $V^\prime \in \VMC^{2,a}_{reg}$. 
Finally $\VMC^a_{reg} = \VMC^{2,a}_{reg} \cap C^\infty (M\times S^1,\R)$ 
implies $V^\prime \in \VMC^a_{reg}$.

  {\bf ad \boldmath $ii)$}
Since $\VMC^a_{reg}$ is open and dense in $(C^\infty , d)$ it is residual.
The identity
\[
  \VMC_{reg} = \bigcap_{a \in \N} \VMC^a_{reg}
\]
implies the claim, because 
any countable intersection of residual sets
is again residual and therefore dense.
\end{proof}

\begin{appendix}
\section{Finite Sum} \label{appsec:finite sum}

    For every $a \in \R$ and every $V \in \VMC^a_{reg}$, the set  $Crit$
  is a $0$-dimensional manifold and the set $Crit^a$ is finite.
  We recall that $Crit^a$ consists by definition of the smooth maps $x:S^1 \to M$ which satisfy
  \begin{equation} \label{eq:*}
    -\n_t \dot x - \n V_t (x) = 0
  \end{equation}
  as well as $\SMC_V(x) <a$. 
  As a consequence the sum in (\ref{eq:sum}) is finite.
  The proof is standard and combines regularity theory, the implicit function theorem
  and compactness arguments.

  \vspace{.2cm}
  \noindent
  {\bf Regularity.}
  We need to extend the domain of definition of $\SMC_V$
  to the Sobolev space $W^{2,2}(S^1,M)$ in order to apply
  the implicit function theorem in the next step.
  Our aim is to show that every $x \in W^{2,2}(S^1,M)$ 
  which satisfies (\ref{eq:*}) almost everywhere
  is necessarily $C^\infty$-smooth.
  In view of the Sobolev embedding theorem we know that every
  $x \in W^{2,2}(S^1,M)$ is indeed of class $C^1$ and $\dot x$
  is absolutely continous.
  If $x$ is in addition a solution to (\ref{eq:*}) we see that in local coordinates
  it holds almost everywhere
  $$
    \ddot x^k = - \Gamma_{ij}^k (x) \dot x^i \dot x^j - g^{k\ell}(x) \frac{\p V_t}{\p x^\ell} (x)
  $$
  Because the right hand side is of class $C^0$ and $\dot x$ is absolutely continous,
  it follows that $x$ is of class $C^2$.
  Hence the right hand side is $C^1$ and so $x$ is $C^3$.
  The iteration continues and we obtain finally $x \in C^\infty$.

  \vspace{.2cm}
  \noindent
  {\bf Implicit function theorem.}
  If $V \in \VMC_{reg}$, we know
  that zero is a regular value of the Fredholm section $\FMC$ in (\ref{eq:F}), which
  is of Fredholm index $0$.
  Hence it follows from the infinite dimensional implicit function theorem that 
  the zero set of $\FMC$
  -- which by regularity is precisely $Crit$ --
  is a submanifold of dimension $0$ of the domain $W^{2,2}(S^1,M)$.
  In particular this means that the elements of $Crit$ are isolated.

  \vspace{.2cm}
  \noindent
  {\bf Compactness.}
  If $V \in \VMC^a_{reg}$, then
  the set  $Crit^a$ consists of finitely many elements: Let us assume by contradiction
  that it contains infinitely many distinct elements $\{ x_\nu \}_{\nu \in \N}$.
  We prove that there exists $x\in Crit^a$ and a subsequence $\{ x_{\nu_k} \}_{k \in \N}$
  converging to $x$ in $W^{2,2}(S^1,M)$, which is a contradiction to the former paragraph.

    As we observed in section \ref{sec:indexthm}, an element $x_\nu \in Crit^a$
  corresponds to a $1$-periodic Hamiltonian orbit
  $z_\nu=z_{x_\nu}=(x_\nu,g(x_\nu) \dot x_\nu)$ and we are going to
  prove that the uniform bound $a$ for $\SMC_V$ implies uniform bounds
  for the initial conditions $(x_\nu^0,y_\nu^0):=(x_\nu(0),g(x_\nu(0)) \dot x_\nu(0)) \in T^*M$.
  So the sequence of initial conditions lies in a compact subset of $T^*M$
  and therefore admits a convergent subsequence
  $(x_{\nu_k}^0, y_{\nu_k}^0) \to (x_0,y_0) \in T^*M$ for $k \to \infty$.
  Let $\varphi_t : T^*M \to T^*M$ be the time-$t$-map of the Hamiltonian flow.
  The $(x_{\nu_k}^0, y_{\nu_k}^0)$ are fixed points
  of $\varphi_1$ and -- because $\varphi_1$ is continous -- so is $(x_0,y_0)$.
  In other words the limit $z(t)=\varphi_t (x_0,y_0)$ is a $1$-periodic orbit.
  Setting $x(t)=\pi (z(t))$, where $\pi : T^*M \to M$ is the natural projection,
  we obtain that $x \in Crit^a$. \\
  This shows that $x_{\nu_k} \to x$ in $C^1$. Moreover, using the fact that $x$ and $x_{\nu_k}$
  both satisfy (\ref{eq:*}) this implies $x_{\nu_k} \to x$ in $C^2$.
  In view of the continous embedding $C^2(S^1,M) \hookrightarrow W^{2,2}(S^1,M)$
  we obtain convergence in $W^{2,2}(S^1,M)$, 
  but this contradicts the fact that the elements of $Crit$ are isolated.
  \\
  It remains to get the uniform bounds: Since $M$ is compact, there is nothing
  to prove for the base components $x_\nu(0) \in M$.
  Now the bound $a$ for the classical action $\SMC_V$ leads to a uniform
  $L^2$-bound for $\dot x$ for all $x \in W^{2,2}(S^1,M)$ with $\SMC_V(x) <a$, namely
  \begin{equation} \label{eq:unif-bound}
    \| \dot x \|_{L^2}^2 
    =\int_0^1 | \dot x |^2 dt < 2a + 2 \| V \|_{L^\infty(S^1 \times M)}
  \end{equation}
  By (\ref{eq:*}), we get
  $\frac{d}{dt} |x_\nu|^2= -2 \langle \n V_t(x_\nu) , \dot x_\nu \rangle$ pointwise in $t$.
  Integrate this identity over the interval $[0,t]$ to obtain
  $$
    | \dot x_\nu (0) |^2 \le | \dot x_\nu (t) |^2
    + \| \n V \|_{L^\infty}^2 \int_0^1 | \dot x_\nu (\tau) |^2 \: d\tau
  $$
  and hence by integrating again and using (\ref{eq:unif-bound}) it follows
  $$
    | \dot x_\nu (0) |^2 \le \left( 1+\| \n V \|_{L^\infty}^2 \right) \| \dot x_\nu \|_{L^2}^2
    < 2 \left( a + \| V \|_{L^\infty} \right) \left( 1+\| \n V \|_{L^\infty}^2 \right)
  $$
  where the right hand side only depends on $a$ and $V$.
\end{appendix}

   


\begin{thebibliography}{99999}

  \bibitem[Br83]{Br83} Brezis H., 
   {\sl Analyse fonctionelle -- Th\`{e}orie et applications}, Masson, Paris 1983. 

  \bibitem[Ci94]{Ci94} Cieliebak K., {\sl Pseudo-holomorphic curves
   and periodic orbits on cotangent bundles},
   J. Math. Pures Appl. {\bf 73} (1994), 251--278.

  \bibitem[CZ84]{CZ84} Conley C., Zehnder E., {\sl Morse-type index theory 
   for flows and periodic solutions for hamiltonian equations}, 
   Comm. Pure Appl. Math. {\bf XXXVII} (1984), 207--253.

  \bibitem[DS94]{DS94} Dostoglou S., Salamon D.A.,
   {\sl Cauchy-Riemann operators, self-duality, and the spectral flow},
   in 'First European Congress of Mathematics' {\bf I}, Invited Lectures (Part 1), 
   Joseph A., Mignot F., Murat F., Prum B., Rentschler R. (editors),
   Birkh\"auser Verlag, Progress in Mathematics {\bf 119} (1994), 511--545.

  \bibitem[Du76]{Du76} Duistermaat J.J., {\sl On the Morse index in variational calculus}, 
   Advances in Math. {\bf 21} (1976), 173--195.

  \bibitem[Fl89]{Fl89} Floer A., {\sl Symplectic fixed points
   and holomorphic spheres} Commun. Math. Phys. {\bf 120} (1989) 575--611.

  \bibitem[GL58]{GL58} Gelfand I.M., Lidskii V.B., {\sl On the structure of the 
   regions of stability of linear canonical systems of differential equations 
   with periodic coefficients}, Translations A.M.S. ({\bf 2}) 8 (1958), 143--181.

  \bibitem[Hi76]{Hi76} Hirsch M.W., {\sl Differential topology},
   Graduate Texts in Math. {\bf 33}, Springer-Verlag New York 1976.

  \bibitem[RS80]{RS80} Reed M., Simon B., {\sl Methods of modern mathematical 
   physics, I Functional analysis}, Academic Press 1980.

  \bibitem[RS93]{RS93} Robbin J., Salamon D.A., {\sl The Maslov index for paths},
   Topology {\bf 32} (1993), 827--844.

  \bibitem[RS95]{RS95} Robbin J., Salamon D.A., {\sl The spectral flow
   and the Maslov index}, Bull. London Math. Soc. {\bf 27} (1995), 1--33.

  \bibitem[Sa96]{Sa96} Salamon D.A., {\sl Spin geometry and Seiberg-Witten 
   invariants}, preliminary version June 1996.

  \bibitem[Sa99]{Sa99} Salamon D.A., {\sl Lectures on Floer Homology},
   in 'Symplectic Geometry and Topology',  Eliashberg and Traynor (editors), 
   IAS/Park City Mathematics series {\bf 7} (1999), 143--230.

  \bibitem[Sm73]{Sm73} Smale S., {\sl An infinite dimensional
   version of Sard's theorem}, Am. J. Math. {\bf 87} (1973), 213--221.

  \bibitem[SW01]{SW01} Salamon D.A., Weber J., {\sl $J$-holomorphic curves
   in cotangent bundles and Morse theory on the loop space}, in preparation.

  \bibitem[Vi90]{Vi90} Viterbo C., {\sl A new obstruction to embedding Lagrangian tori}
   Invent. Math. {\bf 100} (1990), no. 2, 301--320.

  \bibitem[Vi96]{Vi96} Viterbo C., {\sl Functors and computations in Floer
   homology with applications, Part II}, Preprint October 1996, revised February 1998.

  \bibitem[We99]{We99} Weber J., {\sl $J$-holomorphic curves in cotangent
   bundles and the heat flow}, Dissertation TU Berlin, 1999.

  \bibitem[We01]{We01} Weber J., {\sl Geodesic homology}, in preparation.


\end{thebibliography}
\end{document}